\documentclass[12pt]{amsart}
\usepackage[T1]{fontenc}
\usepackage{amsmath}
\usepackage{amsfonts}
\usepackage{amssymb}
\usepackage{xcolor}
\usepackage{bbding}
\usepackage{makecell}
\usepackage{stmaryrd}
\usepackage{txfonts}
\usepackage{graphicx}
\usepackage{diagbox}
\usepackage{epsfig}
\usepackage{xypic}
\usepackage{tikz}
\usepackage{longtable}
\usepackage{caption}
\usepackage{pifont}
\usepackage[all]{xy}
\usepackage{pgflibraryarrows}
\usepackage{pgflibrarysnakes}
\usepackage[shortlabels]{enumitem}
\usepackage{tikz}
\usepackage{tikz-qtree}
\usepackage{ifpdf}
\ifpdf
\usepackage[colorlinks,
linkcolor=magenta,       
anchorcolor=magenta,    
citecolor=magenta,     
final,
hyperindex]{hyperref}
\else
\usepackage[colorlinks,
linkcolor=magenta,       
anchorcolor=magenta,     
citecolor=magenta,    
final,
hyperindex,driverfallback=hypertex]{hyperref}
\fi

\usepackage{xspace}

\usepackage{listings}

\newtheorem{theorem}{Theorem}[section]
\newtheorem{lemma}[theorem]{Lemma}
\newtheorem{coro}[theorem]{Corollary}

\newtheorem{prop}[theorem]{Proposition}
\theoremstyle{definition}

\newtheorem{exmp}[theorem]{Example}

\newcommand{\nc}{\newcommand}
\newcommand{\delete}[1]{}

\nc{\tred}[1]{\textcolor{red}{#1}}
\nc{\tblue}[1]{\textcolor{blue}{#1}} \nc{\tgreen}[1]{\textcolor{green}{#1}} \nc{\tpurple}[1]{\textcolor{purple}{#1}} \nc{\btred}[1]{\textcolor{red}{\bf #1}} \nc{\btblue}[1]{\textcolor{blue}{\bf #1}} \nc{\btgreen}[1]{\textcolor{green}{\bf #1}} \nc{\btpurple}[1]{\textcolor{purple}{\bf #1}}

\newcommand{\efootnote}[1]{}

\nc{\mlabel}[1]{\label{#1}}  
\nc{\mcite}[2][]{\cite[#1]{#2}}  
\nc{\mref}[1]{\ref{#1}}  
\nc{\mbibitem}[1]{\bibitem{#1}} 

\delete{
\nc{\mlabel}[1]{\label{#1}  
{\hfill \hspace{1cm}{\bf{{\ }\hfill(#1)}}}}
\nc{\mcite}[1]{\cite{#1}}  
\nc{\mref}[1]{\ref{#1}{{\bf{{\ }(#1)}}}}  
\nc{\mbibitem}[1]{\bibitem[\bf #1]{#1}} 
}

\renewcommand\geq{\geqslant}
\renewcommand\leq{\leqslant}

\renewcommand\bar[1]{\overline{#1}}

\nc{\nz}{\varepsilon}
\nc{\Id}{\mathrm{Id}}
\nc{\map}[2]{{#2}^{#1}}
\nc{\gp}{B}

\nc{\Irr}{\mathrm{Irr}}
\nc{\vx}{\sigma} \nc{\vy}{\tau} \nc{\dvx}{\sigma^{(1)}} \nc{\dvy}{\tau^{(1)}} \nc{\done}{\vep} \nc{\mcitep}[1]{\mcite{#1}} \nc{\wt}{\mathrm{wt}} \nc{\bre}[1]{|#1|} \nc{\mapmonoid}{\frakM} \nc{\disjoint}{\frakM'}
\nc{\ncpoly}[1]{\langle #1\rangle}  
\nc{\mapm}[1]{\lfloor\!|{#1}|\!\rfloor}
\nc{\diff}[1]{{}^\NC\{ #1 \}} \nc{\disj}[1]{\{{#1}\}'} \nc{\mdisj}[1]{\frakM'(#1)} \nc{\brho}{\bar{\rho}} \nc{\om}{\bar{\frakm}} \nc{\frakn}{\mathfrak n} \nc{\ddeg}[1]{^{(#1)}} \nc{\opset}{X} \nc{\genset}{{Z}} \nc{\NC}{\mathrm{{NC}}} \nc{\leaf}{\mathrm{leaf}} \nc{\twig}{\mathrm{twig}} \nc{\fe}{\mathrm{fl}} \nc{\munderline}[1]{#1} \nc{\bo}{o} \nc{\dep}{\mathrm{depth}} \nc{\ofe}{\mathrm{ofl}} \nc{\dfe}{\mathrm{dfe}} \nc{\fex}{\mathrm{fex}} \nc{\dl}{\mathrm{dlex}} \nc{\db}{\mathrm{db}} \nc{\lex}{\mathrm{lex}} \nc{\clex}{\mathrm{clex}} \nc{\dgp}{\mathrm{dgp}} \nc{\dgx}{\mathrm{dgx}} \nc{\br}{\mathrm{br}} \nc{\obd}{\mathrm{odb}} \nc{\ob}{\mathrm{ob}}
\nc{\pie}{\mathrm{PIE}}
\nc{\rbo}{\mathrm{RBO}}
\nc{\supp}{\mathcal{S}}
\nc{\nul}{\mathcal{Z}}

\nc{\bin}[2]{ (_{\stackrel{\scs{#1}}{\scs{#2}}})}  
\nc{\binc}[2]{ \left (\!\! \begin{array}{c} \scs{#1}\\
    \scs{#2} \end{array}\!\! \right )}  
\nc{\bincc}[2]{  \left ( {\scs{#1} \atop
    \vspace{-1cm}\scs{#2}} \right )}  
\nc{\bs}{\bar{S}} \nc{\cosum}{\sqsubset} \nc{\la}{\longrightarrow} \nc{\rar}{\rightarrow} \nc{\dar}{\downarrow} \nc{\dprod}{**} \nc{\dap}[1]{\downarrow \rlap{$\scriptstyle{#1}$}} \nc{\md}[1]{\bar{#1}} \nc{\uap}[1]{\uparrow \rlap{$\scriptstyle{#1}$}} \nc{\defeq}{\stackrel{\rm def}{=}} \nc{\disp}[1]{\displaystyle{#1}} \nc{\dotcup}{\ \displaystyle{\bigcup^\bullet}\ } \nc{\gzeta}{\bar{\zeta}} \nc{\hcm}{\ \hat{,}\ } \nc{\hts}{\hat{\otimes}} \nc{\barot}{{\otimes}} \nc{\free}[1]{\bar{#1}} \nc{\uni}[1]{\tilde{#1}} \nc{\hcirc}{\hat{\circ}} \nc{\leng}{\ell} \nc{\lleft}{[} \nc{\lright}{]} \nc{\lc}{\lfloor} \nc{\rc}{\rfloor}
\nc{\lb}{[} 
\nc{\rb}{]} 
\nc{\curlyl}{\left \{ \begin{array}{c} {} \\ {} \end{array}
    \right.  \!\!\!\!\!\!\!}
\nc{\curlyr}{ \!\!\!\!\!\!\!
    \left. \begin{array}{c} {} \\ {} \end{array}
    \right \} }
\nc{\longmid}{\left | \begin{array}{c} {} \\ {} \end{array}
    \right. \!\!\!\!\!\!\!}
\nc{\onetree}{\bullet} \nc{\ora}[1]{\stackrel{#1}{\rar}}
\nc{\ola}[1]{\stackrel{#1}{\la}}
\nc{\ot}{\otimes} \nc{\mot}{{{\boxtimes\,}}} \nc{\otm}{\overline{\boxtimes}} \nc{\sprod}{\bullet} \nc{\scs}[1]{\scriptstyle{#1}} \nc{\mrm}[1]{{\rm #1}} \nc{\msum}{\sum\limits}
\nc{\margin}[1]{\marginpar{\rm #1}}   
\nc{\dirlim}{\displaystyle{\lim_{\longrightarrow}}\,} \nc{\invlim}{\displaystyle{\lim_{\longleftarrow}}\,} \nc{\mvp}{\vspace{0.3cm}} \nc{\tk}{^{(k)}} \nc{\tp}{^\prime} \nc{\ttp}{^{\prime\prime}} \nc{\svp}{\vspace{2cm}} \nc{\vp}{\vspace{8cm}} \nc{\proofbegin}{\noindent{\bf Proof: }}
\nc{\proofend}{$\blacksquare$ \vspace{0.3cm}}
\nc{\modg}[1]{\!<\!\!{#1}\!\!>}
\nc{\intg}[1]{F_C(#1)} \nc{\lmodg}{\!<\!\!} \nc{\rmodg}{\!\!>\!} \nc{\cpi}{\widehat{\Pi}}
\nc{\sha}{{\,\makebox[0.6em]{\cyr X}\,}}  
\nc{\shap}{{\mbox{\cyrs X}}} 
\nc{\shpr}{\diamond}    
\nc{\shp}{\ast} \nc{\shplus}{\shpr^+}
\nc{\shprc}{\shpr_c}    
\nc{\msh}{\ast} \nc{\zprod}{m_0} \nc{\oprod}{m_1} \nc{\vep}{\varepsilon} \nc{\labs}{\mid\!} \nc{\rabs}{\!\mid}
\nc{\astarrow}{\overset{\raisebox{-3pt}{$\ast$}}{\rightarrow}}

\nc{\Sym}{{\mrm{Sym}}}
\nc{\Nsym}{{\mrm{NSym}}}
\nc{\Qsym}{{\mrm{QSym}}}
\nc{\SSym}{{\mrm{ SSym}}}
\nc{\SSRCT}{{\mrm{SSRCT}}}
\nc{\SSRCTs}{{\mrm{SSRCTs}}}
\nc{\SRCT}{{\mrm{SRCT}}}
\nc{\SRCTs}{{\mrm{SRCTs}}}
\nc{\SSRT}{{\mrm{SSRT}}}
\nc{\SSRTs}{{\mrm{SSRTs}}}
\nc{\syms}{{symmetric functions\xspace}}
\nc{\qsyms}{{quasi-symmetric functions\xspace}}
\nc{\nsymg}{{\mathrm{NSym}_\gp}}
\nc{\HSym}{{\mrm{\mathfrak{H}Sym}}}
\nc{\parr}{{\mrm {Par}}}
\nc{\Set}{{\mrm {Set}}}
\nc{\Comp}{\mrm {Comp}}
\nc{\comp}{\mrm {comp}}
\nc{\Des}{{\mrm {Des}}}
\nc{\pc}{\mrm {pc}}
\nc{\Sol}{{\mrm {Sol}}}
\nc{\Span}{\mrm {span}}
\nc{\Sh}{{\mrm {Sh}}}
\nc{\st}{{\rm{st}}}

\nc{\dth}{d} \nc{\mmbox}[1]{\mbox{\ #1\ }} \nc{\fp}{\mrm{FP}} \nc{\rchar}{\mrm{char}} \nc{\Fil}{\mrm{Fil}} \nc{\Mor}{Mor\xspace} \nc{\gmzvs}{gMZV\xspace} \nc{\gmzv}{gMZV\xspace} \nc{\mzv}{MZV\xspace} \nc{\mzvs}{MZVs\xspace} \nc{\Hom}{\mrm{Hom}} \nc{\id}{\mrm{id}} \nc{\im}{\mrm{im}} \nc{\incl}{\mrm{incl}}  \nc{\mchar}{\rm char}

\nc{\Alg}{\mathbf{Alg}} \nc{\Bax}{\mathbf{Bax}} \nc{\bff}{\mathbf f} \nc{\bfk}{{\bf k}} \nc{\bfone}{{\bf 1}} \nc{\bfx}{\mathbf x} \nc{\bfy}{\mathbf y}
\nc{\base}[1]{\bfone^{\otimes ({#1}+1)}} 
\nc{\Cat}{\mathbf{Cat}} \delete{}
\nc{\detail}{\marginpar{\bf More detail}
    \noindent{\bf Need more detail!}
    \svp}
\nc{\Int}{\mathbf{Int}} \nc{\Mon}{\mathbf{Mon}}
\nc{\rbtm}{{shuffle }} \nc{\rbto}{{Rota-Baxter }} \nc{\remarks}{\noindent{\bf Remarks: }} \nc{\Rings}{\mathbf{Rings}} \nc{\Sets}{\mathbf{Sets}}
\nc{\balpha}{\mathbf{\alpha}}

\nc{\BA}{{\mathbb A}} \nc{\CC}{{\mathbb C}} \nc{\DD}{{\mathbb D}} \nc{\EE}{{\mathbb E}} \nc{\FF}{{\mathbb F}} \nc{\GG}{{\mathbb G}} \nc{\HH}{{\mathbb H}} \nc{\LL}{{\mathbb L}} \nc{\NN}{{\mathbb N}} \nc{\KK}{{\mathbb K}} \nc{\PP}{{\mathbb P}} \nc{\QQ}{{\mathbb Q}} \nc{\RR}{{\mathbb R}} \nc{\TT}{{\mathbb T}} \nc{\VV}{{\mathbb V}} \nc{\ZZ}{{\mathbb Z}}

\nc{\cala}{{\mathcal A}} \nc{\calc}{{\mathcal C}} \nc{\cald}{{\mathcal D}} \nc{\cale}{{\mathcal E}} \nc{\calf}{{\mathcal F}} \nc{\calg}{{\mathcal G}} \nc{\calh}{{\mathcal H}} \nc{\cali}{{\mathcal I}} \nc{\call}{{\mathcal L}} \nc{\calm}{{\mathcal M}} \nc{\caln}{{\mathcal N}} \nc{\calo}{{\mathcal O}} \nc{\calp}{{\mathcal P}} \nc{\calr}{{\mathcal R}} \nc{\cals}{{\mathcal S}} \nc{\calt}{{\mathcal T}} \nc{\calw}{{\mathcal W}} \nc{\calk}{{\mathcal K}} \nc{\calx}{{\mathcal X}}
\nc{\calz}{{\mathcal Z}}
 \nc{\CA}{\mathcal{A}}

\nc{\fraka}{{\mathfrak a}} \nc{\frakA}{{\mathfrak A}} \nc{\frakb}{{\mathfrak b}} \nc{\frakB}{{\mathfrak B}}
\nc{\frakc}{{\mathfrak c}}  \nc{\frakD}{{\mathfrak D}}
\nc{\frakH}{{\mathfrak H}}
\nc{\frakh}{{\mathfrak h}} \nc{\frakM}{{\mathfrak M}}
\nc{\frakO}{{\mathfrak O}}
\nc{\frakE}{{\mathfrak E}}
\nc{\bfrakM}{\overline{\frakM}} \nc{\frakm}{{\mathfrak m}} \nc{\frakP}{{\mathfrak P}} \nc{\frakN}{{\mathfrak N}} \nc{\frakp}{{\mathfrak p}} \nc{\frakS}{{\mathfrak S}}
\nc{\frakk}{{\mathfrak k}}
\nc{\frakx}{{\mathfrak x}}
\nc{\frakl}{{\mathfrak l}} \nc{\ox}{\bar{\frakx}} \nc{\frakX}{{\mathfrak X}} \nc{\fraky}{{\mathfrak y}} \nc\dop{\delta}
\nc{\Reduce}{{\rm Red}}

\nc{\name}[1]{{\bf #1}}

\font\cyr=wncyr10 \font\cyrs=wncyr7
\nc{\redt}[1]{\textcolor{red}{#1}}
\nc{\yu}[1]{\textcolor{red}{\tt Yu:#1}}
\nc{\liao}[1]{\textcolor{blue}{\tt Liao:#1}}

\begin{document}
\title[Combinatorial equivalence of separable elements]{Combinatorial equivalence of separable elements in types $A$ and $B$}

\author{Yong Liao}
\address{School of Mathematics and Statistics, Southwest University, Chongqing 400715, China}
\email{l2694433238y@email.swu.edu.cn}

\author{Yuping Yang}
\address{School of Mathematics and Statistics, Southwest University, Chongqing 400715, China}
\email{yupingyang@swu.edu.cn}

\author{Houyi Yu$^{\ast}$}  \thanks{*Corresponding author}
\address{School of Mathematics and Statistics, Southwest University, Chongqing 400715, China}
\email{yuhouyi@swu.edu.cn}

\hyphenpenalty=8000

\begin{abstract}
We study the combinatorial equivalence of separable elements in types $A$ and $B$. 
A bijection is constructed from the set of separable permutations in the symmetric group $S_{n+1}$ to  the set of separable signed permutations in the hyperoctahedral group $B_n$. This bijection preserves descent statistics and induces a poset isomorphism under the left weak order. As a consequence, separable signed permutations are enumerated by the large Schr\"oder numbers, and their descent polynomials are shown to be $\gamma$-positive.
Building on a recursive characterization of separable signed permutations via direct sum and skew sum operations, we derive explicit product formulas for the rank generating functions of the principal upper and lower ideals of separable signed permutations under the left weak order. 
\end{abstract}

\keywords{Separable permutation, Signed permutation, Weak order, Rank generating function, Pattern, Descent polynomial}

\maketitle

\hyphenpenalty=8000 \setcounter{section}{0}


\allowdisplaybreaks
\section{Introduction}\label{sec:int}

Let $S_n$ denote the symmetric group of all permutations on the set $[n]=\{1,2,\ldots,n\}$. We represent a permutation $w\in S_n$ in one-line notation as $w=w_1w_2\cdots w_n$.
A permutation $w$ in $S_n$ is called \emph{separable} if it avoids the patterns $3142$ and $2413$. Equivalently, $w$ is separable if there exist no indices $i_1<i_2<i_3<i_4$ such that either $w_{i_2}<w_{i_4}<w_{i_1}<w_{i_3}$ or 
$w_{i_3}<w_{i_1}<w_{i_4}<w_{i_2}$.
Although such permutations were first formally introduced by Bose, Buss, and Lubiw \cite{BBL98} in their study of the permutation pattern matching problem, the concept itself can be traced back to earlier work of Avis and Newborn \cite{AN81} on pop-stacks. 

Owing to their recursive structure, separable permutations possess a wealth of elegant properties.
A key example is their enumeration by the large Schr\"oder numbers~\cite{Sha91,Wes95}. 
Furthermore, the family of separable permutations can be characterized constructively as those permutations that can be generated from the singleton permutation $1$ by recursively applying the operations of \emph{direct sum} and \emph{skew sum} on permutations. For any $u\in S_p$ and $v\in S_q$, the direct sum $u \oplus v$  and the skew sum $u \ominus v$ are defined respectively by
\begin{align*}
	(u \oplus v)_i 
	&=\begin{cases}
		u_i, &\, \text{if } 1\leq i \leq p, \\
		v_{i-p} + p, &\, \text{if } p+1\leq i \leq p+q,
	\end{cases}\\
	(u \ominus v)_i &=
	\begin{cases}
		u_i + q, &\quad \text{if }  1\leq i \leq p,\\
		v_{i-p}, &\quad \text{if }p+1\leq i \leq p+q.
	\end{cases}
\end{align*}
The significance of separable permutations is further highlighted by their frequent appearance in diverse mathematical fields, including bootstrap percolation~\cite{Sha91}, pattern matching~\cite{NRV16}, Schr\"oder paths~\cite{YY24}, probability theory~\cite{BBFGP18},  permutation classes~\cite{AAV11,Vat15}, the distribution of permutation statistics~\cite{CKZ24,FLZ18}, and the splittings of the symmetric group~\cite{BW88,GG20am}.

The concept of separable permutations has recently been extended beyond the symmetric group. Gaetz and Gao \cite{GG20aam} introduced the notion of separable elements for arbitrary Weyl groups, generalizing the defining criteria from the classical case. These elements are characterized by pattern avoidance within the associated root systems, following the framework of Billey and Postnikov~\cite{BP05}. This characterization is profound, as it induces a natural algebraic decomposition of the Weyl group.
A key combinatorial result of Gaetz and Gao \cite{GG20am}  is the product formulas for the rank generating functions of the principal upper and lower order ideals generated by a separable element. These formulas are expressed through the combinatorics of nested sets on the graph associahedron constructed from the Dynkin diagram of the Weyl group.
More recently, Gossow and Yacabi \cite{GY25} demonstrated a deep connection to representation theory by proving that, in simply-laced Weyl groups, the action of separable elements on canonical bases yields a bijection up to lower-order terms.

The purpose of this work is to provide a detailed study of separable elements of type $B$. We focus particularly on their relationship with type $A$ separable elements, and on deriving explicit formulas for the rank generating functions of the principal upper and lower order ideals generated by a type $B$ separable element.
Our results demonstrate a combinatorial equivalence between these two classes of elements.
Our starting point is \cite[Lemma 3.5]{LY25}, which gives an equivalent characterization of type $B$ separable elements in terms of pattern avoidance. 
Recall that the Weyl group of type $B$ can be realized as the hyperoctahedral group, whose elements are called signed permutations. It was shown in
\cite[Lemma 3.5]{LY25} that a signed permutation is separable if and only if it avoids the signed permutation patterns $\bar{2}1$, $2\bar{1}$, $3142$, $2413$, $\bar{3}\,\bar{1}\,\bar{4}\,\bar{2}$, and $\bar{2}\,\bar{4}\,\bar{1}\,\bar{3}$. 

Our first main results establish several equivalent characterizations and fundamental properties of separable signed permutations (Theorem \ref{lem:separablesignedperm}). In particular, we generalize the direct sum and skew sum operations to signed permutations and prove that a signed permutation is separable if and only if it can be constructed from the base signed permutations $1$ and $\bar{1}$ by recursively applying these operations.

Building on this characterization, we construct a bijection $\varphi_n$ from the set $K(S_{n+1})$ of separable permutations in $S_{n+1}$ to the set $K(B_n)$ of separable signed permutations in $B_n$ (Theorem \ref{thm:bijectionKSn+1KBn}). This bijection immediately implies that the cardinality of $K(B_n)$ is given by the $n$th Schr\"oder number.
The bijection $\varphi_n$ possesses two key properties that highlight a certain equivalence between separable permutations and separable signed permutations.
First, $\varphi_n$ preserves descents and double descents. As a consequence, and building on a result of Fu, Lin, and Zeng \cite[Theorem 1.1]{FLZ18}, we prove that
$$
\sum_{w \in K( B_n)} t^{\mathrm{des}_B(w)}=\sum_{w \in K(S_{n+1})} t^{\text{des}_A(w)},
$$
a polynomial that is $\gamma$-positive, and hence palindromic and unimodal (Theorem \ref{thm:descentpolyB}). Second,  $\varphi_n$ is order-preserving with respect to the left weak order, establishing a poset isomorphism between $K(S_{n+1})$ and $K(B_n)$ under this order (Theorem \ref{thm:order-preserving}).

Finally, we show that for any separable signed permutation $w$, the rank generating functions of the principal upper and lower ideals $\Lambda_{w}^B$ and $V_{w}^B$ in the left weak order can be reduced to those of its standardized permutation counterpart. Combined with a result of Wei \cite[Theorem 3.5]{We12}, 
this leads to explicit product formulas for the rank generating functions of $\Lambda_{w}^B$ and $V_{w}^B$.

The paper is organized as follows. Section \ref{sec:preliminaries} reviews necessary background on hyperoctahedral groups and symmetric groups.
Section \ref{sec:Separablesignedpermutations} presents equivalent characterizations of separable signed permutations. In Section \ref{Sec:relationSn+1Bn}, we construct the bijection between $K(S_{n+1})$ and $K(B_{n})$ and study its properties.
Finally, Section \ref{sec:formulasforFlambdaVB} is devoted to the derivation of explicit formulas for rank generating functions of the left weak order ideals.

\section{Preliminaries}\label{sec:preliminaries}

This section recalls some basic facts about hyperoctahedral groups and symmetric groups; additional details can be found in \cite{BB05, Hum90, Sta12}.

Unless otherwise specified, $n$ denotes a nonnegative integer.
For integers $m$ and $n$ with $m\leq n$, we write $[m,n]=\{m,m+1,\ldots,n\}$, and use $[n]$ as shorthand for the set $[1,n]$.
For convenience the additive inverse of an integer $n$
is often denoted by $\bar{n}$.

Let $(W, S)$ be a finite Coxeter system. Thus, $W$ is a group generated by $S$, subject only to relations of the form $(ss')^{m_{ss'}}=1$, where $m_{ss}=1$
and $m_{ss'}=m_{s's}\geq 2$ for $s\neq s'$ in $S$. The group $W$ is a \emph{Weyl group} if $m_{ss'}\in\{2,3,4,6\}$ when $s\neq s'$.
Each element $w$ in $W$ can be written as a product   $w=s_{i_1}s_{i_2}\cdots s_{i_r}$ with $s_{i_j}\in S$.
If $r$ is minimal, then such an expression is called a \emph{reduced expression}, and $r$ is called the \emph{length} of $w$, denoted $\ell(w)=r$.

Let $P$ be a poset.  For each element $x\in P$, the \emph{principal upper order ideal} generated by $ x $ is $V_x = \{ y \in P \mid y \geq x \}$, and the \emph{principal lower order ideal} generated by $ x $ is  $ \Lambda_x = \{ y \in P \mid y \leq x \} $. If $ P $ is a graded poset of rank $n$ with $p_i$ elements of rank $i$, then the polynomial
\begin{align*}
	F(P, q) = \sum_{i=0}^{n} p_i q^i
\end{align*}
is called the \emph{rank generating function} of $P$.

The \emph{left weak order} on a Weyl group $W$ is defined by $u \leq_L w$ if $\ell(w)=\ell(u)+\ell(wu^{-1})$. This order makes the group $W$ a graded poset ranked by the length, and in fact a lattice, with minimal element $e$ and maximal element $w_0$. Thus, for any $w \in W$, we have $V_ {w}=[w, w_0]$ and $\Lambda_ {w}=[e, w]$. If $W$ has rank $n$, then its rank generating function \cite{Hum90} is given by 
\begin{align*}
	F(W, q) = \sum_{w \in W} q^{\ell(w)} = \prod_{i=1}^{n} \frac{q^{d_i} - 1}{q - 1},
\end{align*}
where the $d_i$ are the \emph{degrees} of $W$.

Throughout this paper, we focus on Weyl groups of types $A$ and $B$.
Following \cite[Chapters 1 and 8]{BB05}, the Weyl group of type $A_{n-1}$ is isomorphic to the symmetric group $S_n$ of all permutations of $[n]$. Analogously,
the Weyl group of type $B_n$ is isomorphic to the \emph{hyperoctahedral group} $B_n$, which consists of all \emph{signed permutations} on $[n]$, where a signed permutation is a bijection $w$ on the set
$[\pm1,\ldots,\pm n]$ such that $w(\bar i)=\bar {w(i)}$ for all $i\in[n]$. We identify $S_n$ as a subgroup of $B_n$ in the natural way.

A signed permutation may be represented in one-line notation
$w=w_1w_2\cdots w_n$, where $w_i=w(i)$ for $i\in[n]$, since the images of the negative entries in $w$ are then uniquely determined. Let $\varepsilon$ denote the empty (signed) permutation, and let $ w_0^{S_n}$ and  $w_0^{B_n} $ be the longest elements of $ S_n$  and $ B_n$, respectively. Then $ w_0^{ S_n}=n(n-1)\cdots 1$ and  $w_0^{B_n}=\overline{1}\,\overline{2}\cdots\overline{n}$.
For a signed permutation $w$, we denote by $\vert w \vert$ the permutation $\vert w_1\vert \vert w_2 \vert \cdots \vert w_n \vert$, and by $\bar{w}$ the signed permutation $\bar{w_1}\,\bar{w_2} \cdots \bar{w_n}$. Note that $\overline{w}=w_0^{B_n}w=ww_0^{ B_n}$. For example, if $w=3\overline{5}\,\overline{2}41\overline{6}$, then
$|w|=352416$ and $\overline{w}=\overline{3}52\overline{4}\,\overline{1}6$.

Given a word $a=a_1a_2\cdots a_n$ over the integers, its \emph{standardized permutation}, denoted $\mathrm{st}(a)$, is the permutation $w$ in $S_n$ with the same relative order as $a$. That is,
for all $i,j$ with $1\leq i<j\leq n$, we have $w_{i}<w_{j}$ if and only if $a_i\leq a_j$.
If all $a_i\neq0$, the \emph{standardized signed permutation} of $a$,
denoted $\mathrm{sts}(a)$, is the signed permutation $w\in B_n$ such that $\mathrm{Neg}(w)=\mathrm{Neg}(a)$, and for all $i,j$ with $1\leq i<j\leq n$,
$|w_i|<|w_j|$ if and only if $|a_i|\leq|a_j|$, where $\mathrm{Neg}(a)=\{i\in[n]\mid a_i<0\}$.
For example, $\mathrm{st}(4\overline{7}\,\overline{3}52\overline{9})=523641$,
$\mathrm{sts}(4\overline{7}\,\overline{3}52\overline{9})=3\overline{5}\,\overline{2}41\overline{6}$.
Note that if all entries $w_i$ are positive integers, then $\mathrm{st}(w)=\mathrm{sts}(w)$.

Let $m,n$ be positive integers with $m\leq n$, and let $u\in  B_m$, $w\in  B_n$.
We say $w$ \emph{contains the pattern} $u$ if there exist
$1\leq i_1<i_2<\cdots<i_m\leq n$ such that 
$\mathrm{sts}(w_{i_1}w_{i_2}\cdots w_{i_m})=u$.
Otherwise, $w$ is said to \emph{avoid the pattern} $u$.
For example, consider the signed permutation $w=\bar{2}53\bar{6}\,\bar{4}1$.
Taking the subsequence $w_3w_4w_6=3\bar{6}1$, we have $\mathrm{sts}(w_3w_4w_6)=2\bar{3}1$, so $w$ contains $2\bar{3}1$ as a pattern. But $w$ avoids $3\bar{1}2$.

The study of signed patterns arises in various contexts, including fully commutative signed permutations \cite{Ste97}, characterizations of rationally smooth Schubert varieties \cite{Bil98}, and vexillary signed permutations whose Stanley symmetric functions are single Schur $Q$-functions \cite{BL98}. More information on this topic can be found in \cite{Kit11}.

\section{Separable signed permutations}\label{sec:Separablesignedpermutations}

Separable elements in general Weyl groups were originally defined by Gaetz and Gao \cite{GG20aam} and were further studied in \cite{GG20am}.
In the specific case of type $B$, these elements were characterized combinatorially in \cite[Lemma 3.5]{LY25}, which shows that a signed permutation
 $w\in B_n$ is \emph{separable}
if and only if it avoids the patterns $\bar{2}1$, $2\bar{1}$, $3142$, $2413$, $\bar{3}\,\bar{1}\,\bar{4}\,\bar{2}$, and $\bar{2}\,\bar{4}\,\bar{1}\,\bar{3}$. 
Consequently, a permutation $w\in S_n$ is separable in the classical sense if and only if it remains separable when viewed as an element of $ B_n$.
Denoting by $K(S_n)$ and $K(B_n)$ the sets of separable elements in $S_n$ and $B_n$, respectively,
it then follows that $K(S_n)= K(B_n)\cap S_n$.

\begin{lemma}[\cite{GG20am}, Corollary 2]\label{lem:ww-1rl}
	Let $W$ be a Weyl group. If $w\in W$ is separable, then $w_0w$, $ww_0$, and
	$w^{-1}$ are also separable. 
\end{lemma}

\begin{lemma}\label{lem:|w|}
 For any integer $n\geq1$, $K( S_n)=\{|w|\mid w\in K( B_n)\}$.
\end{lemma}
	\begin{proof}
	The inclusion $K( S_n)\subseteq \{|w|\mid w\in K( B_n)\}$ is immediate since $K( S_n)\subseteq  K( B_n)$ and $w=|w|$ for all $w\in S_n$.
		
	For the reverse inclusion, let $w\in K( B_n)$, and let $u=|w|$, so $u\in S_n$.
		Suppose on the contrary that $u\notin K( S_n)$. Then $u$ contains either the pattern $2413$ or $3142$. Since $2413$ is the group inverse of $3142$, by Lemma \ref{lem:ww-1rl},
		 we can suppose without loss of generality that $u$ contains the patterns $2413$.
		Thus, there exist indices $i_1<i_2<i_3<i_4$ such that $u_{i_3}<u_{i_1} <u_{i_4}<u_{i_2}$, implying $|w_{i_3}|<|w_{i_1}| <|w_{i_4}|<|w_{i_2}|$.
		Since $w$ avoids $\overline{2}1$ and $ 2\overline{1}$,
		the subsequence $w_{i_1}w_{i_2}w_{i_3}w_{i_4}$ also avoids $\overline{2}1$ and $ 2\overline{1}$.
		 
		If $w_{i_1}>0$, then avoiding $2\bar{1}$ forces $w_{i_3}>0$, and avoiding $\bar{2}1$ forces $w_{i_2}>0$, which in turn implies $w_{i_4}>0$ since $w$ avoids $2\bar{1}$.
		But then $\mathrm{sts}(w_{i_1}w_{i_2}w_{i_3}w_{i_4})=2413$, contradicting the separability of $w$. 		
		If $w_{i_1}<0$, consider $\bar{w}=w_0^{B_n}w$, which is separable by Lemma \ref{lem:ww-1rl}. A completely analogous argument shows that $\bar{w}$ contains $2413$, again a contradiction. Thus, $u\in K(S_n)$, and hence
		 $K( S_n)=\{|w|\mid w\in K( B_n)\}$.
	\end{proof}
	 
This result indicates that separable signed permutations in $B_n$ arise from separable permutations in $ S_n$  by assigning appropriate signs to their entries.

We now introduce two key operations that recursively construct separable signed permutations.
For $u\in B_p$ and $v\in B_q$, define the \emph{direct sum}  $u\oplus v\in B_{p+q}$ by
\begin{align*}
	(u \oplus v)_i 
	=&\begin{cases}
		u_i, & \text{if } 1\leq i \leq p, \\
		v_{i-p} + (\mathrm{sgn}\, v_{i-p})p, & \text{if } p+1 \leq i \leq  p+q,
	\end{cases}
\end{align*}
where  $\mathrm{sgn}\,x=1$ if $x>0$, and $-1$ if $x<0$.
The \emph{skew sum} $u\ominus v$ is defined only when all entries of $u$ and $v$ have the same sign (i.e. $\mathrm{sgn}\,u_i=\mathrm{sgn}\,v_j$ for all $i\in[p]$ and $j\in[q]$). In that case,
\begin{align*}
	(u \ominus v)_i =
	\begin{cases}
		u_i + (\mathrm{sgn\,} u_i)q, & \text{if } 1\leq i \leq p,\\
		v_{i-p}, & \text{if }  p+1 \leq i \leq  p+q.
	\end{cases}
\end{align*}
For example,  $$1\bar{2}\,\bar{3}\oplus\bar{2}\,\bar{1}=1\bar{2}\,\bar{3}\,\bar{5}\,\bar{4}\quad \text{and}\quad \bar{1}\,\bar{2}\,\bar{3}\ominus\bar{2}\,\bar{1}=\bar{3}\,\bar{4}\,\bar{5}\,\bar{2}\,\bar{1}.$$ However, the expression $1\bar{2}\,\bar{3}\ominus\bar{2}\,\bar{1}$ is undefined due to inconsistent signs.

The following classical result \cite{BBJJS11} characterizes separable permutations constructively.
\begin{lemma}\label{lem:seppermutopluominu}
	A permutation is separable if and only if it can be constructed from the permutation $1$ by repeated application of the operations $\oplus$ and $\ominus$.
\end{lemma}

A recursive decomposition property of separable permutations follows directly.

\begin{lemma}\label{lem:separableproperty}
	A permutation $w\in S_n$ is separable if and only if there exist $p\in[n-1]$ and separable permutations $u\in S_p$, $v\in S_{n-p}$ such that either
	$w=u\oplus v$ or $w=u\ominus v$.
\end{lemma} 

Recall from \cite{AS05} that a permutation $w\in S_n$  has a \emph{global ascent} (respectively, \emph{global descent}) at position $p\in[n-1]$ if $w_i<w_j$ (respectively, $w_i>w_j$) for all $i\in[p]$ and $j\in [p+1,n]$.
Let $\mathrm{GAsc}(w)$ and $\mathrm{GDes}(w)$ denote the sets of global ascents and global descents of $w$, respectively.
For example, 
$$\mathrm{GAsc}(3125467)=\{3,5,6\} \quad\text{and}\quad
\mathrm{GDes}(5746312)=\{4,5\}.$$ 
Note that $p\in \mathrm{GAsc}(w)$ if and only if $w=w_1\cdots w_p\oplus \mathrm{st}(w_{p+1}\cdots w_n)$, and $p\in \mathrm{GDes}(w)$ if and only if $w=\mathrm{st}(w_1\cdots w_p)\ominus w_{p+1}\cdots w_n$.
By Lemma \ref{lem:separableproperty},
every separable permutation of length greater than $1$ contains at least one global ascent or a global descent.

The following lemma provides a decomposition for separable signed permutations.

\begin{lemma}\label{lem:septypeBopluominu2}
	Let $w\in K(B_n)$ and $p\in[n-1]$ such that $w_p$ and $w_j$ have opposite signs for every $j\in[p+1,n]$. Then $w_1\cdots w_p$ and $\mathrm{sts}(w_{p+1}\cdots w_n)$ are separable with $$w=w_1\cdots w_p\oplus \mathrm{sts}(w_{p+1}\cdots w_n).$$
\end{lemma}
\begin{proof}
	The separability of the factors is clear since they are patterns of $w$. Assume first that $w_p<0$, so $w_{j}>0 $  for $j>p$. Since $w$ avoids $\bar{2}1$ and $2\bar{1}$, 
for any $i\in[p]$ we have  
$|w_i|<w_j$ for all $j>p$. Thus, the direct sum decomposition holds. If $w_p>0$, then $w_{j}<0$  for all $j>p$. By Lemma \ref{lem:ww-1rl}, $\bar{w}=w_0^{ B_n}w$ is a separable signed permutation.
Applying the preceding proof for the case where $w_p<0$ to $\bar{w}$, yields that
$$\bar{w}=\bar{w}_1\cdots \bar{w}_p\oplus \mathrm{sts}(\bar{w}_{p+1}\cdots \bar{w}_n),$$
and hence 
$$w=w_1\cdots w_p\oplus \mathrm{sts}(w_{p+1}\cdots w_n),$$
completing the proof.
\end{proof}

We now present the main theorem of this section, characterizing separable signed permutations in multiple ways.

\begin{theorem}\label{lem:separablesignedperm}
	For $w\in  B_n$, the following statements are equivalent:
	\begin{enumerate}
		\item\label{item:wsepsignedp1} $w$ is a separable signed permutation.
		\item\label{item:wsepsignedp2op}  There exist $p\in[n-1]$, separable signed permutations $u\in  B_p$, $v\in  B_{n-p}$ such that either
		$w=u\oplus v$ or $w=u\ominus v$.
		\item\label{item:wsepsignedstsop} There exists a unique sequence $0=p_0< p_1<p_2<\cdots <p_k< p_{k+1}=n$ such that
		\begin{align*}
		w=w_1\cdots w_{p_1}\oplus \mathrm{sts}(w_{p_1+1}\cdots w_{p_2})\oplus\cdots\oplus \mathrm{sts}(w_{p_{k}+1}\cdots w_{n}),
		\end{align*}
where $w_1\cdots w_{p_1}$ and $\mathrm{sts}(w_{p_{i}+1}\cdots w_{p_{i+1}})$  for $i\in [k]$ are separable, and for each $i\in[k]$, the consecutive blocks $w_{p_{i-1}+1}\cdots w_{p_i}$ and $w_{p_{i}+1}\cdots w_{p_{i+1}}$ have entries of opposite signs.
		\item\label{item:wsepsignedp1bar1} $w$ can be built from $1$ and $\bar1$ by repeated applications of the operations $\oplus $ and $\ominus $.
	\end{enumerate}
\end{theorem} 
\begin{proof}
	\eqref{item:wsepsignedp1} $\Rightarrow$ \eqref{item:wsepsignedp2op}
If
$\mathrm{Neg}(w)=\emptyset$, then $w\in K(S_n)$ and the result follows from Lemma \ref{lem:separableproperty}. If $\mathrm{Neg}(w)\subsetneqq[n]$ is nonempty, then \eqref{item:wsepsignedp2op} follows from Lemma \ref{lem:septypeBopluominu2}.  If $\mathrm{Neg}(w)=[n]$, then by Lemma \ref{lem:ww-1rl},  $w_{0}^{ B_{n}}w$ is a separable permutation, and hence there exist separable  permutations $u\in  S_p$ and $v\in  S_{n-p}$ for some $p\in[n]$ such that either $w_0^{ B_{n}}w=u\oplus v$ or $w_0^{ B_{n}}w=u\ominus v$.
Write $u=u_1\cdots u_p$ and $v=v_1\cdots v_{n-p}$, then
\begin{align*}
	w_0^{ B_{n}}(u \oplus v)
	&=w_0^{ B_{n}}(u_1\cdots u_p (v_{1}+p)\cdots(v_{n-p}+p))\\
	&=\overline{u_1}\cdots \overline{u_p}\, \overline{v_{1}+p}\cdots\overline{v_{n-p}+p}\\
	&=w_0^{ B_{p}}u \oplus w_0^{ B_{n-p}}v,	
\end{align*}
and similarly $w_0^{ B_{n}}(u \ominus v)=w_0^{ B_{p}}u \ominus w_0^{ B_{n-p}}v.$ Thus,  either 
$$w=w_0^{ B_{n}}(u\oplus v)=w_0^{ B_{p}}u \oplus w_0^{ B_{n-p}}v\quad \text{or}\quad 
w=w_0^{ B_{n}}(u\ominus v)=w_0^{ B_{p}}u \ominus w_0^{ B_{n-p}}v,$$
with both components separable by Lemma \ref{lem:ww-1rl}.
Hence, \eqref{item:wsepsignedp2op} follows.

	 \eqref{item:wsepsignedp2op} $\Rightarrow$ \eqref{item:wsepsignedstsop}
	 The proof is by induction on $n$. If $n=1$, then $w$ is either $1$ or $\bar{1}$, and the direct sum decomposition is trivial. For $n\geq2$, by Item \eqref{item:wsepsignedp2op}, we may assume that $w=u\oplus v$ or $w=u\ominus v$, where $u\in  K(B_p)$ and $v\in K(B_{n-p})$ for some $p\in[n-1]$. Apply the induction hypothesis to $u$ and $v$ to get their unique block decomposition, say,
	 \begin{align}
	 	u&=u_1\cdots u_{p_1}\oplus \mathrm{sts}(u_{p_1+1}\cdots u_{p_2})\oplus\cdots\oplus \mathrm{sts}(u_{p_{k}+1}\cdots u_{p}),\label{eq:udirectsum1}\\
	 	v&=v_1\cdots v_{q_1}\oplus \mathrm{sts}(v_{q_1+1}\cdots v_{q_2})\oplus\cdots\oplus \mathrm{sts}(v_{q_{k}+1}\cdots v_{n-p}).\label{eq:udirectsum2}
	 \end{align}
	 
	 For the case $w=u\oplus v$, if $u_p$ and $v_1$ have the same sign, then 
	 \begin{align*}
	 	w=u_1\cdots u_{p_1}\oplus \mathrm{sts}(u_{p_1+1}\cdots u_{p_2})\oplus\cdots&\oplus \left(\mathrm{sts}(u_{p_{k}+1}\cdots u_{p})\oplus v_1\cdots v_{q_1}\right)\\
	 	&\oplus \mathrm{sts}(v_{q_1+1}\cdots v_{q_2})\oplus\cdots\oplus \mathrm{sts}(v_{q_{k}+1}\cdots v_{n-p});
	 \end{align*}
	 if  $u_p$ and $v_1$ have different signs, then $w$ is the direct sum of the right hand sides of Eqs.\,\eqref{eq:udirectsum1} and \eqref{eq:udirectsum2}.
	 
	 For the case $w=u\ominus v$, the entries of $u$ and $v$ must have the same signs. Hence $\mathrm{Neg}(w)$ is either $\emptyset$ or $[n]$, and taking $p_1=n$ is enough. This proves \eqref{item:wsepsignedstsop}.
 
\eqref{item:wsepsignedstsop} $\Rightarrow$ \eqref{item:wsepsignedp1bar1} By induction on $n$.
	 Since $1$ and $\bar{1}$ are separable, \eqref{item:wsepsignedp1bar1} is true for $n=1$. Assume that $n\geq2$. If $\mathrm{Neg}(w)=\emptyset$, then $w$ is a separable permutation, so, by Lemma \ref{lem:seppermutopluominu}, $w$ can be built from the permutation $1$ by applying $ \oplus $ and $ \ominus $ repeatedly.
	 If $\mathrm{Neg}(w)=[n]$, then $\overline{w}$ is a separable permutation, so by Lemma \ref{lem:seppermutopluominu} again, $\overline{w}$ can be built from the permutation $1$, and hence $w$ can be built from the permutation $\overline{1}$, by applying $ \oplus $ and $ \ominus $ repeatedly. If $\mathrm{Neg}(w)$ is neither $\emptyset$ nor $[n]$, then  the proof follows straightforwardly from \eqref{item:wsepsignedstsop} by induction on $n$.
	 
	 \eqref{item:wsepsignedp1bar1}  $\Rightarrow$  \eqref{item:wsepsignedp1} 
	If $w$ is built from $1$ and $\bar1$ via $ \oplus $ and $ \ominus$, 
	then $|w|$ is a separable permutation built from $1$, so $|w|$ avoids $3142$ and $2413$, and hence $w$ avoids $3142$, $2413$, $\bar{3}\,\bar{1}\,\bar{4}\,\bar{2}$, and $\bar{2}\,\bar{4}\,\bar{1}\,\bar{3}$. Since for any $u\in S_p$ and $v\in S_{q}$, the skew sum $u\ominus v$ is only defined when the entries $u_i$ and $v_j$ have the same sign for all $i\in[p]$ and $j\in[q]$, we see that $w$ avoids the patterns $\bar{2}1$ and $2\bar{1}$. Thus, $w$ is a separable signed permutation, completing the proof.
\end{proof}

 The (skew) direct sum decompositions of a separable signed permutation are not necessarily unique in general. For example, $1243=1\oplus 132=12\oplus 21$, and
 $12\bar3\,\bar4=12\oplus \bar{1}\,\bar{2}=12 \bar{3}\oplus\bar{1}$.
 However, the decomposition in Theorem \ref{lem:separablesignedperm}\eqref{item:wsepsignedstsop} is unique; for the example above, the decompositions are $1243$ and $12\oplus \bar{1}\,\bar{2}$, respectively.

We conclude this section with a useful relation between separable signed permutations and their standardizations, which will be utilized in subsequent arguments.
\begin{coro}\label{lem:st(w)inksn}
	If $w \in K(B_n)$, then $\mathrm{st}(w) \in K(S_n)$.
\end{coro}
\begin{proof}
	Let $w=w_1w_2\cdots w_n\in K(B_n)$. First consider the case where $w_n<0$.
	If $\mathrm{Neg}(w)=[n]$, then all entries $w_i<0$. Since $w$ avoids $\bar{3}\,\bar{1}\,\bar{4}\,\bar{2}$ and $\bar{2}\,\bar{4}\,\bar{1}\,\bar{3}$, we see that $\mathrm{st}(w)$ avoids $2413$ and $3142$, and thus $\mathrm{st}(w)\in K(S_n)$.
	If  $\mathrm{Neg}(w)\neq[n]$, then, by Lemma \ref{lem:septypeBopluominu2}, there exists $p\in[n]$ such that
	$w=w_1\cdots w_p \oplus \mathrm{sts}(w_{p+1}\cdots w_n)$, where $w_{i}<0$ for all $i\in[p+1,n]$. Thus,
	$$\mathrm{st}(w)=\mathrm{st}(w_1\cdots w_p) \ominus \mathrm{st}(w_{p+1}\cdots w_n).$$
	Since $w \in K(B_n)$, it follows that $w_1\cdots w_p$ and $w_{p+1}\cdots w_n$ are separable signed permutations, so that $\mathrm{st}(w_1\cdots w_p)$ and $\mathrm{st}(w_{p+1}\cdots w_n)$ are separable permutations by induction.
	Hence $\mathrm{st}(w)\in K(S_n)$ by Lemma \ref{lem:separableproperty}.
	For the case where $w_n>0$, we have $$\mathrm{st}(w)=w_0^S\mathrm{st}(w_0^Bw).$$ 
	Since $w_0^Bw\in K(B_n)$ with last entry negative, we see that $\mathrm{st}(w_0^Bw)\in K(S_n)$ by the previous case and hence  $\mathrm{st}(w)\in K(S_n)$ by Lemma \ref{lem:ww-1rl}.
\end{proof}

\section{Combinatorial equivalence between $K( S_{n+1})$ and $K( B_n)$}\label{Sec:relationSn+1Bn}
In this section, we construct a bijection from $K( S_{n+1})$ to $K( B_n)$ that is also an isomorphism under the left weak order,
thereby establishing a combinatorial equivalence between these two classes of separable elements.
This bijection implies that separable signed permutations are enumerated by the large Schr\"oder numbers. Furthermore, it preserves the descent and double descent statistics, ensuring the equality of the descent polynomials for $K( S_{n+1})$ and $K( B_n)$. This common polynomial is  $\gamma$-positive, and hence palindromic and unimodal.

\subsection{A bijection between $K( S_{n+1})$ and $K( B_n)$}\label{Sec:bijection}

For each nonnegative integer $n$, we define recursively  a bijection  $\varphi_n:K( S_{n+1})\rightarrow K( B_n)$ as follows.
Let the empty signed permutation $\varepsilon\in  B_0$ act as the identity element for both direct sum and skew sum operations.
Define $\varphi_0(1)=\varepsilon$, and for $n\geq1$, set 
\begin{align}\label{eq:defnvarphin}
	\varphi_n(w)=\begin{cases}
		\varphi_{p-1}(w_1\cdots w_p)\oplus\mathrm{st}(w_{p+1}\cdots w_{n+1}), & \mathrm{if}\  p\in\mathrm{GAsc}(w),\\
		\varphi_{p-1}(\mathrm{st}(w_1\cdots w_p))\oplus w_0^{ B_{n+1-p}} w_0^{ S_{n+1-p}}(w_{p+1}\cdots w_{n+1}),& \mathrm{if}\  p\in\mathrm{GDes}(w).
	\end{cases}
\end{align}

For example, $\varphi_1(12)=1$, $\varphi_1(21)=\bar1$. The values of $\varphi_n(w)$ for $n=2$ and $3$ are listed in Table $1$ and Table $2$, respectively.

\vspace{5mm}

	\begin{center}
		\ \ Table $1$: The values of $\varphi_2(w)$.\\
		\vspace{3mm}
		\renewcommand{\arraystretch}{1.8} 
		\setlength{\tabcolsep}{6pt}  
		\begin{tabular}{*{7}c}
			\hline
			$w$&$123$&$132$&$213$&$231$&$312$&$321$\\
			$\varphi_2(w)$&$12$&$21$&$\bar12$&$1\bar2$&$\bar2\,\bar1$&$\bar1\,\bar2$\\
			\hline
		\end{tabular}
	\end{center}
	\vspace{5mm}
	
		\begin{center}
		\ \ Table $2$: The values of $\varphi_3(w)$.\\
		\vspace{3mm}
		 \renewcommand{\arraystretch}{1.8} 
		 \setlength{\tabcolsep}{6pt}  
		\begin{tabular}{*{12}c}
			\hline
			$w$&$1234$&$2134$&$1324$&$1243$&$3124$&$2314$&$2143$&$1423$&$1342$&$4123$&$3214$\\
			$\varphi_3(w)$&$123$&$\bar123$& $213$&$132$&$\bar2\,\bar13$&$1\bar23$&$\bar132$&$312$&$231$&$\bar3\,\bar2\,\bar1$&$\bar1\,\bar23$\\
			\hline
			$w$&$2341$&$1432$&$4213$&$4132$&$3412$&$3241$&$2431$&$4312$&$4231$&$3421$&$4321$\\
			$\varphi_3(w)$& $12\bar3$&$321$&$\bar2\,\bar3\,\bar1$&$\bar3\,\bar1\,\bar2$&$1\bar3\,\bar2$&$\bar12\bar3$&$21\bar3$&$\bar1\,\bar3\,\bar2$&$\bar2\,\bar1\,\bar3$&$1\bar2\,\bar3$&$\bar1\,\bar2\,\bar3$\\
			\hline
		\end{tabular}
	\end{center}
\vspace{5mm}
	
Note that global ascents and global descents of a permutation are not necessarily unique. For example, in the permutation $4231$, both positions $1$ and $3$ are global descents. Nevertheless, the following lemma shows that the value of $\varphi_n(w)$ does not depend on the particular choice of global ascents (or global descents) in $w$, ensuring that $\varphi_n$ is well-defined.
	
\begin{lemma}\label{lem:varphinwelldefine}
For any nonnegative integer $n$, the map $\varphi_n$ is well-defined.
\end{lemma}
\begin{proof}
We proceed by induction on $n$. The base case $n=0$ is trivial. Now assume that $\varphi_i$ is well-defined for all $i<n$.
	Let $w\in S_{n+1}$, and suppose $p,q\in[n]$ with $p<q$ are both global ascents (or global descents) of $w$. By the induction hypothesis, both $\varphi_{p-1}$ and $\varphi_{q-1}$ are well-defined.
	
	If $\mathrm{GAsc}(w)$ is nonempty, we assume that $p,q\in\mathrm{GAsc}(w)$, then  
	\begin{align*}
	\{w_1,\ldots,w_p\}=[p],\quad \{w_{p+1},\ldots,w_q\}=[p+1,q],\quad \{w_{q+1},\ldots,w_{n+1}\}=[q+1,n+1],
\end{align*}
which implies $p\in \mathrm{GAsc}(w_1\cdots w_q)$ and
\begin{align*}
	\mathrm{st}(w_{p+1}\cdots w_{n+1}) =\mathrm{st}(w_{p+1}\cdots w_{q})\oplus\mathrm{st}(w_{q+1}\cdots w_{n+1}).
\end{align*} 
Therefore,
	\begin{align*}
		\varphi_{q-1}(w_1\cdots w_q)\oplus\mathrm{st}(w_{q+1}\cdots w_{n+1}) 
		&=\varphi_{p-1}(w_1\cdots w_p)\oplus\mathrm{st}(w_{p+1}\cdots w_{q})\oplus\mathrm{st}(w_{q+1}\cdots w_{n+1}) \\
		&=\varphi_{p-1}(w_1\cdots w_p)\oplus\mathrm{st}(w_{p+1}\cdots w_{n+1}).
	\end{align*}
	
	If $\mathrm{GDes}(w)$ is nonempty, we assume that $p,q\in \mathrm{GDes}(w)$, then 
	\begin{align*}
		\{w_1,\ldots,w_p\}=[n+2-p,n+1],\quad
		\{w_{p+1},\ldots,w_q\}=[n+2-q,n+1-p]
	\end{align*}
	and 
		\begin{align*}
		\quad \{w_{q+1},\ldots,w_{n+1}\}=[n+1-q],
	\end{align*} 
	so that $p\in \mathrm{GDes}(\mathrm{st}(w_{1}\cdots w_{q}))$. Consequently,
	\begin{align*}
		&\,\varphi_{q-1}(\mathrm{st}(w_1\cdots w_q))\oplus w_0^{ B_{n+1-q}} w_0^{ S_{n+1-q}}(w_{q+1}\cdots w_{n+1})\\
		=&\,\varphi_{p-1}(\mathrm{st}(w_1\cdots w_p))\oplus w_0^{ B_{q-p}} w_0^{ S_{q-p}}(\mathrm{st}(w_{p+1}\cdots w_{q}))\oplus w_0^{ B_{n+1-q}} w_0^{ S_{n+1-q}}(w_{q+1}\cdots w_{n+1})\\
		=&\,\varphi_{p-1}(\mathrm{st}(w_1\cdots w_p))\oplus w_0^{ B_{n+1-p}} w_0^{ S_{n+1-p}}(w_{p+1}\cdots w_{n+1}).
	\end{align*}
	
	In both cases, the value of $\varphi_n(w)$ is independent of the choice of global ascents or global descents. Moreover, by Lemma \ref{lem:ww-1rl} and Theorem \ref{lem:separablesignedperm}, we have $\varphi_n(w)\in K(B_n)$. Therefore, $\varphi_n$ is well-defined.
\end{proof}

We now recursively define a map  $\psi_n:K( B_{n})\rightarrow K( S_{n+1})$ and show that it is the inverse of $\varphi_n$, whence it follows that $\varphi_n$ is a bijection.  

Given $w\in K(B_n)$, where $n\geq0$, by Theorem \ref{lem:separablesignedperm}, we can take a direct sum decomposition $w=w_1\cdots w_{p}\oplus \mathrm{sts}(w_{p+1}\cdots w_{n})$, where  $p\in [0,n-1]$, $w_1\cdots w_{p}\in K( B_{p})$,  $\mathrm{sts}(w_{p+1}\cdots w_{n})\in K( B_{n-p})$, and $w_{p+1},\ldots, w_{n}$ have the same sign.
Note that when $p=0$, the direct sum decomposition of $w$ is trivial, that is, 
$w=\varepsilon\oplus w$.
Let $\psi_0(\varepsilon)=1$, and for $n\geq1$, define $\psi_n:K( B_{n})\rightarrow K( S_{n+1})$ by 
\begin{align}\label{defn:psi_n}
	\psi_n(w)=\begin{cases}
		\psi_{p}(w_1\cdots w_{p})\oplus \mathrm{sts}(w_{p+1}\cdots w_{n})& \mathrm{if}\ w_n>0,\\
		\psi_{p}(w_1\cdots w_{p})\ominus w_0^{ S_{n-p}}w_0^{ B_{n-p}}\mathrm{sts}(w_{p+1}\cdots w_{n})& \mathrm{if}\ w_n<0.
	\end{cases}
\end{align}
The following lemma guarantees that $\psi_n(w)$ is independent of the choice of the direct sum decompositions of $w$, so the mapping is well-defined.
For example, since
$$312\overline{4}\,\overline{6}\,\overline{5}=312\oplus\overline{1}\,\overline{3}\,\overline{2}=312\overline{4}\oplus\overline{2}\,\overline{1},$$
it follows that 
\begin{align*}
	\psi_6(312\overline{4}\,\overline{6}\,\overline{5})=\psi_3(312)\ominus w_0^{ S_{3}}w_0^{ B_{3}}\overline{1}\,\overline{3}\,\overline{2}
	=(\psi_0(\varepsilon)\oplus 312)\ominus 312=4756312,
\end{align*}
or equivalently,
\begin{align*}
	\psi_6(312\overline{4}\,\overline{6}\,\overline{5})=\psi_4(312\overline{4})\ominus w_0^{ S_{2}}w_0^{ B_{2}}\overline{2}\,\overline{1}
	=(\psi_3(312)\ominus w_0^{ S_{1}}w_0^{ B_{1}}\overline{1})\ominus12=4756312.
\end{align*}

\begin{lemma}\label{lem:psinwelldefine}
	For any nonnegative integer $n$, the map $\psi_n$ is well-defined.
\end{lemma}
\begin{proof}
	Induction on $n$, with the case $n=0$ being true by the definition of $\psi_0$.
	Suppose the result is true for all $i<n$.
	If the choice of $p$ in Eq.\,\eqref{defn:psi_n} is unique, then $\psi_n(w)$ is uniquely determined by induction.
    If the choice of $p$ in Eq.\,\eqref{defn:psi_n} is not unique, then suppose without loss of generality that
	there exists $q>p$ in $[0,n-1]$ such that $w$ has direct sum decompositions
	\begin{align*}
		w=w_1\cdots w_{p}\oplus \mathrm{sts}(w_{p+1}\cdots w_{n})=w_1\cdots w_{q}\oplus \mathrm{sts}(w_{q+1}\cdots w_{n}).
	\end{align*}
	Then
\begin{align*}
	\{w_1,\ldots,w_p\}\subseteq[\overline{p},\overline{1}]\cup[1,p],\quad 
	\{w_{p+1},\ldots,w_q\}\subseteq [\overline{q},\overline{p+1}]\cup [p+1,q],
\end{align*} 
	and
	\begin{align*}
		\{w_{q+1},\ldots,w_n\}\subseteq [\overline{n},\overline{p+1}]\cup[q+1,n],
	\end{align*}
where $w_{p+1},\ldots,w_n$ have the same sign.
	Moreover, $\psi_p$ and $\psi_q$ are well-defined by the induction hypothesis.
	If $w_n$ and hence $w_q$ are positive integers, then 
	\begin{align*}
	\psi_{q}(w_1\cdots w_{q})\oplus \mathrm{sts}(w_{q+1}\cdots w_{n})
		&=\psi_{p}(w_1\cdots w_{p})\oplus \mathrm{sts}(w_{p+1}\cdots w_{q})\oplus \mathrm{sts}(w_{q+1}\cdots w_{n})\\
		&=	\psi_{p}(w_1\cdots w_{p})\oplus \mathrm{sts}(w_{p+1}\cdots w_{n}).
	\end{align*}
	If $w_n$ and hence $w_q$ are negative integers, then 
	\begin{align*}
		&\,\psi_{q}(w_1\cdots w_{q})\ominus  w_0^{ S_{n-q}}w_0^{ B_{n-q}}\mathrm{sts}(w_{q+1}\cdots w_{n})\\
		=&\, \left(\psi_{p}(w_1\cdots w_{p})\ominus  w_0^{ S_{q-p}}w_0^{ B_{q-p}}\mathrm{sts}(w_{p+1}\cdots w_{q})\right)\ominus  w_0^{ S_{n-q}}w_0^{ B_{n-q}}\mathrm{sts}(w_{q+1}\cdots w_{n})\\
		=&\,\psi_{p}(w_1\cdots w_{p})\ominus w_0^{ S_{n-p}}w_0^{ B_{n-p}}\mathrm{sts}(w_{p+1}\cdots w_{n}).
	\end{align*}
Therefore, $\psi_n(w)\in K(S_{n+1})$ is independent of the choice of $p$ and hence $\psi_n$ is well-defined.
\end{proof}
	
\begin{theorem}\label{thm:bijectionKSn+1KBn}
	For any nonnegative integer $n$, the maps $\varphi_n$ and  $\psi_n$ are bijections that are inverses of each other.
\end{theorem}
\begin{proof}
We show that $\varphi_n\psi_n$ and $\psi_n\varphi_n$ are identity maps by induction on $n$.
The case $n=0$ is clear. Assume the result holds for all indices less than $n$.
Let $w=w_1\cdots w_{n+1}$ be an element of $K(S_{n+1})$.
First suppose $p\in \mathrm{GAsc}(w)$. 
It follows from Eqs.\,\eqref{eq:defnvarphin} and \eqref{defn:psi_n} that 
\begin{align*}
	\psi_n\varphi_n(w)&=\psi_n\left(\varphi_{p-1}(w_1\cdots w_p)\oplus\mathrm{st}(w_{p+1}\cdots w_{n+1})\right)\\
	&=\psi_{p-1}\varphi_{p-1}(w_1\cdots w_p)\oplus\mathrm{sts}(\mathrm{sts}(w_{p+1}\cdots w_{n+1})).
\end{align*}
By the induction hypothesis, $\psi_{p-1}\varphi_{p-1}$ is the identity map, so we have
\begin{align*}
	\psi_n\varphi_n(w)=w_1\cdots w_p\oplus\mathrm{sts}(w_{p+1}\cdots w_{n+1})
	=w.
\end{align*}
Similarly, if $p\in \mathrm{GDes}(w)$, then	
\begin{align*}
	\psi_n\varphi_n(w)
	&=\psi_n\left(\varphi_{p-1}(\mathrm{st}(w_1\cdots w_p))\oplus w_0^{ B_{n+1-p}} w_0^{ S_{n+1-p}}(w_{p+1}\cdots w_{n+1})\right)\\
	&=\psi_{p-1}\varphi_{p-1}(\mathrm{st}(w_1\cdots w_p))\ominus w_0^{ S_{n+1-p}}w_0^{ B_{n+1-p}}( w_0^{ B_{n+1-p}} w_0^{ S_{n+1-p}}(w_{p+1}\cdots w_{n+1}))\\
	&=\mathrm{st}(w_1\cdots w_p)\ominus w_{p+1}\cdots w_{n+1}\\
	&=w.
\end{align*}	
Thus, $\psi_n\varphi_n$ is the identity map on the set  $K( S_{n+1})$.

Now assume that $w$ is an element of $K( B_{n})$. By Theorem \ref{lem:separablesignedperm}, we can write $w=w_1\cdots w_{p}\oplus \mathrm{sts}(w_{p+1}\cdots w_{n})$ for some  $p\in [0,n-1]$, where $w_1\cdots w_{p}\in K( B_{p})$,  $\mathrm{sts}(w_{p+1}\cdots w_{n})\in K( B_{n-p})$, and $w_{p+1},\ldots, w_{n}$ have the same sign.
It follows from Eq.\,\eqref{defn:psi_n} that $p\in \mathrm{GAsc}(\psi_n(w))$  if $w_n>0$, and that $p\in \mathrm{GDes}(\psi_n(w))$  if $w_n<0$. So from Eq.\,\eqref{eq:defnvarphin} and the induction hypothesis we obtain that
\begin{align*}
\varphi_n\psi_n(w)
	&=\begin{cases}
		\varphi_p\psi_{p}(w_1\cdots w_p)\oplus\mathrm{st}(\mathrm{sts}(w_{p+1}\cdots w_{n})), & \mathrm{if}\ w_n>0,\\
		\varphi_p\psi_{p}(w_1\cdots w_p)\oplus w_0^{ B_{n-p}} w_0^{S_{n-p}} w_0^{S_{n-p}} w_0^{ B_{n-p}}\mathrm{sts}(w_{p+1}\cdots w_{n})),& \mathrm{if}\  w_n<0
	\end{cases}\\
	&=w_1\cdots w_p \oplus \mathrm{sts}(w_{p+1}\cdots w_{n})\\
	&=w.
\end{align*}
Consequently, $\varphi_n\psi_n$ is the identity map on the set  $K( B_{n})$, completing the proof.
\end{proof}

The number of separable permutations in $ S_n$
is the $(n-1)$th \emph{Schr\"oder number}, which is the sequence \cite[A006318]{Sloane} beginning with $1, 2, 6, 22, 90, 394, 1806$. This result was conjectured by Shapiro and Getu, and proved by West \cite{Wes95} in 1995. Several other proofs of this enumeration formula can be found in \cite{AAV11,Kre00,Sta94}.
As a direct consequence of Theorem \ref{thm:bijectionKSn+1KBn}, we obtain the following corollary.

\begin{coro}
	The set $K( B_n)$ is counted by the $n$th Schr\"oder number.
\end{coro}
\begin{proof}
	It follows from Theorem \ref{thm:bijectionKSn+1KBn} that $\varphi_n$ is a bijection from $K( S_{n+1})$ to $K( B_n)$. Therefore, $\#K(B_n)=\#K(S_{n+1})$, which is the $n$th Schr\"oder number.
\end{proof}

\subsection{Descent polynomials for separable signed permutations}
Let $w=w_1w_2\cdots w_n$ be an element of $ B_n $. An index $i\in[0,n]$ is called a \emph{B-descent} of $w$ if 
$w_i>w_{i+1}$, and a \emph{double B-descent} if
$w_{i-1}>w_i>w_{i+1}$, where $w_0=0$ and $w_{\bar 1}=w_{n+1}=+\infty$.
Denote by $\mathrm{Des}_B(w)$ and $\mathrm{DDes}_B(w)$ the sets of B-descents and double B-descents of $w$, respectively.
Let $\mathrm{des}_B(w)=\#\mathrm{Des}_B(w)$,  $\mathrm{dd}_B(w)=\#\mathrm{DDes}_B(w)$.
For example, if $w=312\overline{4}\,\overline{6}\,\overline{5}$, then
$\mathrm{Des}_B(w)=\{1,3,4\}$, $\mathrm{DDes}_B(w)=\{4\}$, and hence $\mathrm{des}_B(w)=3$,  $\mathrm{dd}_B(w)=1$.

If $w$ is an element of $ S_n$, then an index $i\in[n]$ is called an \emph{A-descent} of $w$ if $w_i>w_{i+1}$, and a \emph{double A-descent} if
$w_{i-1}>w_i>w_{i+1}$, where $w_0=w_{n+1}=+\infty$.
So if $1$ is an A-descent of a permutation $w$, then it must be a double A-descent of $w$, but not a double B-descent of $w$. We have the similar notation 
$\mathrm{Des}_A(w)$, $\mathrm{DDes}_A(w)$, $\mathrm{des}_A(w)$,  $\mathrm{dd}_A(w)$ for the permutation $w$.

\begin{lemma}\label{prop:w=st0variw}
	For any $w \in K( S_{n+1})$, we have $w=\mathrm{st}(0\cdot \varphi_n(w))$, where $\cdot$ indicates the concatenation of words.
\end{lemma}
	\begin{proof}
	 	We proceed by induction on $n$, where the base case $n=0$ is clear.
	 	Suppose $w=\mathrm{st}(0\cdot \varphi_n(w))$ for all $w \in K( S_{i+1})$,
	 	where $i<n$. Now assume that $w \in K( S_{n+1})$.
	 	If $w=w_1\cdots w_{p}\oplus \mathrm{st}(w_{p+1}\cdots w_{n+1})$ for some $p\in[n-1]$, then $p\in \mathrm{GAsc}(w)$.
	 	It follows from Eq.\,\eqref{eq:defnvarphin} that
	 	\begin{align*}
	 	\varphi_n(w)=\varphi_{p-1}(w_1\cdots w_p)\oplus\mathrm{st}(w_{p+1}\cdots w_{n+1}).
	 	\end{align*}
	 	Since $w_i>0$ for all $i\in[n]$, we get
	 	\begin{align*}
	 		\mathrm{st}(0\cdot \varphi_n(w))&=\mathrm{st}(0\cdot \varphi_{p-1}(w_1\cdots w_p))\oplus\mathrm{st}(w_{p+1}\cdots w_{n+1})\\
	 		&=\mathrm{st}(w_1\cdots w_p)\oplus\mathrm{st}(w_{p+1}\cdots w_{n+1})\\
	 		&=w,
	 	\end{align*}
	 	where the second equality follows from the induction hypothesis.
	 	Similarly, if $w=\mathrm{st}(w_1\cdots w_p)\ominus w_{p+1}\cdots w_{n+1}$, then $p\in \mathrm{GDes}(w)$.
	 	By Eq.\,\eqref{eq:defnvarphin}, we get
	 	\begin{align*}
	 		\varphi_n(w)=\varphi_{p-1}(\mathrm{st}(w_1\cdots w_p))\oplus w_0^{ B_{n+1-p}} w_0^{ S_{n+1-p}}(w_{p+1}\cdots w_{n+1}).
	 	\end{align*}
	 	Since 
	 	\begin{align*}
	 	w_0^{ B_{n+1-p}} w_0^{ S_{n+1-p}}(w_{p+1}\cdots w_{n+1})
	 	=(\overline{n+2-p-w_{p+1}})\cdots (\overline{n+2-p-w_{n+1}}),
	 	\end{align*}
	 	it follows that 
	 	\begin{align*}
	 		\mathrm{st}(0\cdot \varphi_n(w))&=\mathrm{st}(0\cdot \varphi_{p-1}(\mathrm{st}(w_1\cdots w_p)))\ominus\mathrm{st}(w_{p+1}\cdots w_{n+1})\\
	 		&=\mathrm{st}(w_1\cdots w_p)\ominus\mathrm{st}(w_{p+1}\cdots w_{n+1})\\
	 		&=w,
	 	\end{align*}
	 	completing the proof.
	\end{proof}
	For example, $\varphi_5(562341)=1\bar4\,\bar3\,\bar2\,\bar5$, and we have $\mathrm{st}(01\bar4\,\bar3\,\bar2\,\bar5)=562341$.

Alternatively, Lemma \ref{prop:w=st0variw} is equivalent to the following statement since 
$\psi_n$ is the inverse of $\varphi_n$ by Theorem \ref{thm:bijectionKSn+1KBn}.
\begin{coro}\label{cor:psi_nw=st0w}
For any $w \in K( B_{n})$, we have $\psi_n(w)=\mathrm{st}(0\cdot w)$, where $\cdot$ indicates the concatenation of words.
\end{coro}

The following corollary is an immediate consequence of Lemma \ref{prop:w=st0variw} or Corollary
\ref{cor:psi_nw=st0w}.

\begin{lemma}\label{cor:DdesAB}
	For any $w\in K(S_{n+1})$, we have
	\begin{align*}
		\mathrm{Des}_A(w)=\{i+1\mid i\in \mathrm{Des}_B(\varphi_n(w)\},\qquad
		\mathrm{DDes}_A(w)=\{i+1\mid i\in \mathrm{DDes}_B(\varphi_n(w)\}.
	\end{align*}
\end{lemma}
\begin{proof}
	Let $u=\varphi_n(w)$.
 It follows from Lemma \ref{prop:w=st0variw} that $w=\mathrm{st}(0\cdot u)$, so 
 $u_i>u_{i+1}$ if and only if $w_{i+1}>w_{i+2}$, giving $\mathrm{Des}_A(w)=\{i+1\mid i\in \mathrm{Des}_B(\varphi_n(w)\}$. The second identity is proved analogously.
\end{proof}

Recall that the \emph{descent polynomial for separable permutations} in $S_{n}$ is defined as 
$$
	 S_n(t) = \sum_{w \in K(S_n)} t^{\text{des}_A(w)}.
$$
Fu, Lin, and Zeng \cite[Theorem 1.1]{FLZ18} proved that the polynomial $ S_n(t)$ is $\gamma$-positive, and hence is palindromic and unimodal.
More precisely,
\begin{align*}
	 S_n(t) =\sum_{k = 0}^{\left\lfloor \frac{n-1}{2} \right\rfloor} \gamma_{n,k}^{ S}  t^k (1 + t)^{n-1-2k},
\end{align*}
where $\gamma_{n,k}^{ S} = \# \{ w \in K( S_n) \mid \mathrm{dd}_A(w) = 0, \mathrm{des}_A(w) = k \}$, and
$\left\lfloor \frac{n-1}{2} \right\rfloor$ is the greatest integer less than
or equal to $\frac{n-1}{2}$.
An alternative proof was provided by Yang and Yang \cite{YY24} via enumerating the pseudo Schr\"oder paths.

Define the \emph{descent polynomial for separable signed permutations} in $B_n$ as
$$
 B_n(t) = \sum_{w \in K( B_n)} t^{\mathrm{des}_B(w)}.
$$
According to Theorem \ref{thm:bijectionKSn+1KBn} and Lemma \ref{cor:DdesAB}, we have the following result for the polynomial $ B_n(t)$.

\begin{theorem}\label{thm:descentpolyB} 
	We have $B_n(t) =S_{n+1}(t)$, and consequently, 
	\begin{align*}
		 B_n(t) =\sum_{k = 0}^{\left\lfloor \frac{n}{2} \right\rfloor} \gamma_{n,k}^{ B}  t^k (1 + t)^{n-2k},
	\end{align*}
	where $\gamma_{n,k}^{ B} = \# \{ w \in K( B_n) \mid \text{dd}_B(w) = 0, \text{des}_B(w) = k \}=\gamma_{n+1,k}^{ S}.$
	Therefore, the polynomial $ B_n(t)$ is $\gamma$-positive and hence is palindromic and unimodal.
\end{theorem}
\begin{proof}
By Theorem \ref{thm:bijectionKSn+1KBn},
$\varphi_n$ is a bijection from $K( S_{n+1})$ to $K( B_n)$, which together with
Lemma \ref{cor:DdesAB} implies that $\mathrm{des}_A(w)=\mathrm{des}_B(\varphi_n{w})$ for any $w\in K( S_{n+1})$.
Thus, $ B_n(t)= S_{n+1}(t)$ 
and 
$$\# \{ w \in K( B_n) \mid \text{dd}_B(w) = 0, \text{des}_B(w) = k \} = \# \{ w \in K( S_{n+1}) \mid \text{dd}_A(w) = 0, \text{des}_A(w) = k \},$$
so that $\gamma_{n,k}^{ B}=\gamma_{n+1,k}^{ S}$, completing the proof.	
\end{proof}
 
\subsection{Poset isomorphism under the left weak order}
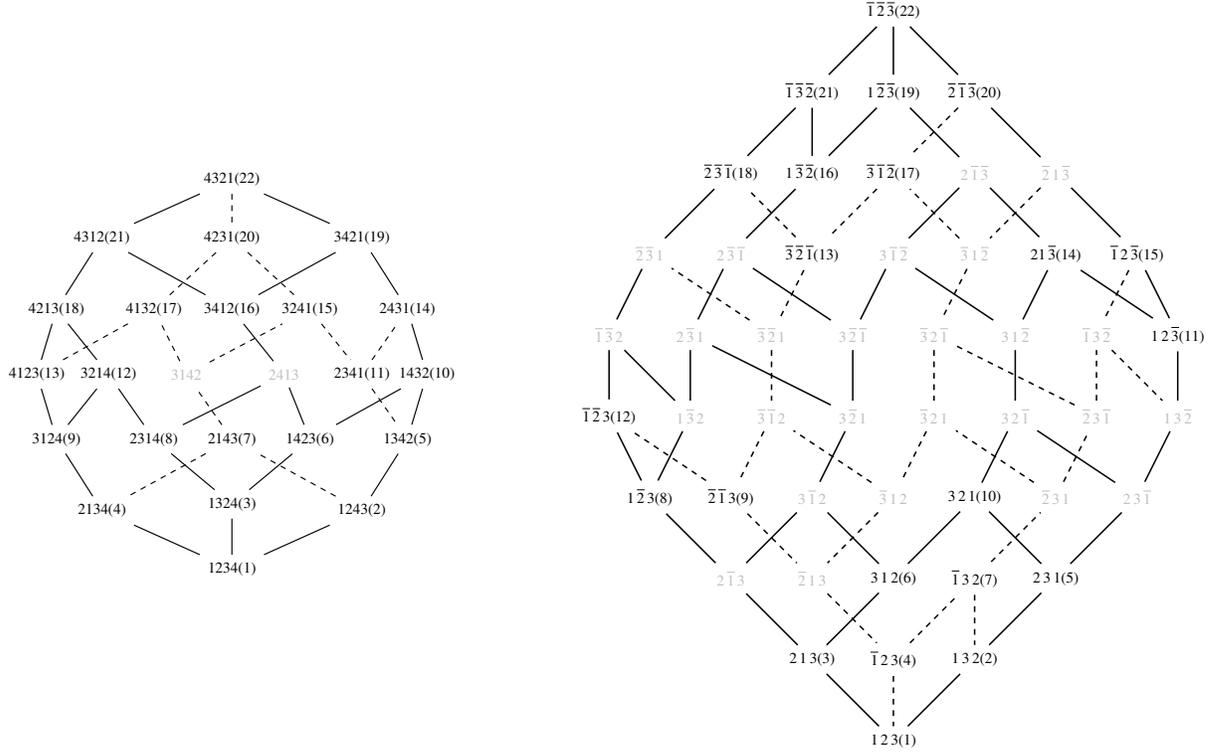
\begin{figure}[ht]
	\centering
	\begin{minipage}[c]{0.45\textwidth}
		\centering
		\scalebox{0.72}{
			\begin{tikzpicture}
				\node (4321) at (0,3.6) {{\tiny$4321(22)$}};
				\node (4312) at (-2.4,2.52) {{\tiny$4312(21)$}};   \node (4231) at (0,2.52) {{\tiny$4231(20)$}}; \node (3421) at (2.4,2.52) {{\tiny$3421(19)$}};
				\node (4213) at (-3.24,1.2) {{\tiny$4213(18)$}};\node (4132) at (-1.44,1.2) {{\tiny$4132(17)$}}; \node (3412) at (0,1.2) {{\tiny$3412(16)$}};
				\node (3241) at (1.44,1.2) {{\tiny$3241(15)$}};  \node (2431) at (3.24,1.2) {{\tiny$2431(14)$}};
				\node (4123) at (-3.6,0) {{\tiny$4123(13)$}};\node (3214) at (-2.28,0) {{\tiny$3214(12)$}};\node (3142) at (-0.84,0) {{\tiny\color{lightgray}$3142$}};
				\node (2413) at (0.96,0) {{\tiny\color{lightgray}$2413$}};
				\node (2341) at (2.4,0) {{\tiny$2341(11)$}};\node (1432) at (3.6,0) {{\tiny$1432(10)$}};
				\node (3124) at (-3.24,-1.2) {{\tiny$3124(9)$}};\node (2314) at (-1.44,-1.2) {{\tiny$2314(8)$}};\node (2143) at (0,-1.2) {{\tiny$2143(7)$}};
				\node (1423) at (1.44,-1.2) {{\tiny$1423(6)$}};\node (1342) at (3.24,-1.2) {{\tiny$1342(5)$}};
				\node (2134) at (-2.4,-2.52) {{\tiny$2134(4)$}};\node (1324) at (0,-2.4) {{\tiny$1324(3)$}};\node (1243) at (2.4,-2.52) {{\tiny$1243(2)$}};
				\node (1234) at (0,-3.6) {{\tiny$1234(1)$}};
				\draw[solid,line width=0.6pt](4321)--(4312)--(4213)--(4123)--(3124)--(2134)--(1234)--(1243)--(1342)--(1432)--(2431)--(3421)--(4321);
				\draw[solid,line width=0.6pt](4312)--(3412)--(2413)--(1423)--(1324)--(1234)
				(1324)--(2314)--(3214)--(4213)        (3124)--(3214)       (2314)--(2413)       (1423)--(1432) (3412)--(3421);
				\draw[dashed,line width=0.6pt](4321)--(4231)--(3241)--(2341)--(1342)   (2341)--(2431)
				(4132)--(4231)
				(3142)--(3241)
				(1243)--(2143)--(3142)--(4132)--(4123)
				(2143)--(2134);
			\end{tikzpicture}
		}
	\end{minipage}
	\hfill
	\begin{minipage}[c]{0.50\textwidth}
		\centering
		\scalebox{0.72}{
			\begin{tikzpicture}
				\node (bar{1}bar{2}bar{3}) at (0,6.75) {{\tiny$\bar{1}\thinspace\bar{2}\thinspace \bar{3}(22)$}};
				\node (bar{1}bar{3}bar{2}) at (-1.5,5.25) {{\tiny$\bar{1}\thinspace\bar{3}\thinspace \bar{2}(21)$}};
				\node (1bar{2}bar{3}) at (0,5.25) {{\tiny$1\thinspace\bar{2}\thinspace \bar{3}(19)$}};
				\node (bar{2}bar{1}bar{3}) at (1.5,5.25) {{\tiny$\bar{2}\thinspace\bar{1}\thinspace \bar{3}(20)$}};
				\node (bar{2}bar{3}bar{1}) at (-3,3.75) {{\tiny$\bar{2}\thinspace\bar{3}\thinspace \bar{1}(18)$}};
				\node (1bar{3}bar{2}) at (-1.5,3.75) {{\tiny$1\thinspace\bar{3}\thinspace \bar{2}(16)$}};
				\node (bar{3}bar{1}bar{2}) at (0,3.75) {{\tiny$\bar{3}\thinspace\bar{1}\thinspace \bar{2}(17)$}};
				\node (2bar{1}bar{3}) at (1.5,3.75) {{\tiny\color{lightgray}$2\thinspace\bar{1}\thinspace \bar{3}$}};
				\node (bar{2}1bar{3}) at (3,3.75) {{\tiny\color{lightgray}$\bar{2}\thinspace1\thinspace \bar{3}$}};
				\node (bar{2}bar{3}1) at (-4.5,2.25) {{\tiny\color{lightgray}$\bar{2}\thinspace\bar{3}\thinspace1$}};
				\node (2bar{3}bar{1}) at (-3,2.25) {{\tiny\color{lightgray}$2\thinspace\bar{3}\thinspace\bar{1}$}};
				\node (bar{3}bar{2}bar{1}) at (-1.5,2.25) {{\tiny$\bar{3}\thinspace\bar{2}\thinspace\bar{1}(13)$}};
				\node (3bar{1}bar{2}) at (0,2.25) {{\tiny\color{lightgray}$3\thinspace\bar{1}\thinspace\bar{2}$}};
				\node (bar{3}1bar{2}) at (1.5,2.25) {{\tiny\color{lightgray}$\bar{3}\thinspace1\thinspace\bar{2}$}};
				\node (21bar{3}) at (3,2.25) {{\tiny$21\thinspace\bar{3}(14)$}};
				\node (bar{1}2bar{3}) at (4.5,2.25) {{\tiny$\bar{1}\thinspace2\thinspace\bar{3}(15)$}};
				\node (bar{1}bar{3}2) at (-5.25,0.75) {{\tiny\color{lightgray}$\bar{1}\thinspace\bar{3}\thinspace2$}};
				\node (2bar{3}1) at (-3.75,0.75) {{\tiny\color{lightgray}$2\thinspace\bar{3}\thinspace1$}};
				\node (bar{3}bar{2}1) at (-2.25,0.75) {{\tiny\color{lightgray}$\bar{3}\thinspace\bar{2}\thinspace1$}};
				\node (3bar{2}bar{1}) at (-0.75,0.75) {{\tiny\color{lightgray}$3\thinspace\bar{2}\thinspace\bar{1}$}};
				\node (bar{3}2bar{1}) at (0.75,0.75) {{\tiny\color{lightgray}$\bar{3}\thinspace2\thinspace\bar{1}$}};
				\node (31bar{2}) at (2.25,0.75) {{\tiny\color{lightgray}$3\thinspace1\thinspace\bar{2}$}};
				\node (bar{1}3bar{2}) at (3.75,0.75) {{\tiny\color{lightgray}$\bar{1}\thinspace3\thinspace\bar{2}$}};
				\node (12bar{3}) at (5.25,0.75) {{\tiny$1\thinspace2\thinspace\bar{3}(11)$}};
				\node (bar{1}bar{2}3) at (-5.25,-0.75) {{\tiny$\bar{1}\thinspace\bar{2}\thinspace3(12)$}};
				\node (1bar{3}2) at (-3.75,-0.75) {{\tiny\color{lightgray}$1\thinspace\bar{3}\thinspace2$}};
				\node (bar{3}bar{1}2) at (-2.25,-0.75) {{\tiny\color{lightgray}$\bar{3}\thinspace\bar{1}\thinspace2$}};
				\node (3bar{2}1) at (-0.75,-0.75) {{\tiny\color{lightgray}$3\thinspace\bar{2}\thinspace1$}};
				\node (bar{3}21) at (0.75,-0.75) {{\tiny\color{lightgray}$\bar{3}\thinspace2\thinspace1$}};
				\node (32bar{1}) at (2.25,-0.75) {{\tiny\color{lightgray}$3\thinspace2\thinspace\bar{1}$}};
				\node (bar{2}3bar{1}) at (3.75,-0.75) {{\tiny\color{lightgray}$\bar{2}\thinspace3\thinspace\bar{1}$}};
				\node (13bar{2}) at (5.25,-0.75) {{\tiny\color{lightgray}$1\thinspace3\thinspace\bar{2}$}};
				\node (1bar{2}3) at (-4.5,-2.25) {{\tiny$1\thinspace\bar{2}\thinspace3(8)$}};
				\node (bar{2}bar{1}3) at (-3,-2.25) {{\tiny$\bar{2}\thinspace\bar{1}\thinspace3(9)$}};
				\node (3bar{1}2) at (-1.5,-2.25) {{\tiny\color{lightgray}$3\thinspace\bar{1}\thinspace2$}};
				\node (bar{3}12) at (0,-2.25) {{\tiny\color{lightgray}$\bar{3}\thinspace1\thinspace2$}};
				\node (321) at (1.5,-2.25) {{\tiny$3\thinspace2\thinspace1(10)$}};
				\node (bar{2}31) at (3,-2.25) {{\tiny\color{lightgray}$\bar{2}\thinspace3\thinspace1$}};
				\node (23bar{1}) at (4.5,-2.25) {{\tiny\color{lightgray}$2\thinspace3\thinspace\bar{1}$}};
				\node (2bar{1}3) at (-3,-3.75) {{\tiny\color{lightgray}$2\thinspace\bar{1}\thinspace3$}};
				\node (bar{2}13) at (-1.5,-3.75) {{\tiny\color{lightgray}$\bar{2}\thinspace1\thinspace3$}};
				\node (312) at (0,-3.75) {{\tiny$3\thinspace1\thinspace2(6)$}};
				\node (bar{1}32) at (1.5,-3.75) {{\tiny$\bar{1}\thinspace3\thinspace2(7)$}};
				\node (231) at (3,-3.75) {{\tiny$2\thinspace3\thinspace1(5)$}};
				\node (213) at (-1.5,-5.25) {{\tiny$2\thinspace1\thinspace3(3)$}};
				\node (bar{1}23) at (0,-5.25) {{\tiny$\bar{1}\thinspace2\thinspace3(4)$}};
				\node (132) at (1.5,-5.25) {{\tiny$1\thinspace3\thinspace2$(2)}};
				\node (123) at (0,-6.75) {{\tiny$1\thinspace2\thinspace3(1)$}};
				\draw[solid,line width=0.8pt] (bar{1}bar{2}bar{3})--(bar{1}bar{3}bar{2})--(bar{2}bar{3}bar{1})--(bar{2}bar{3}1)--(bar{1}bar{3}2)--(bar{1}bar{2}3)--(1bar{2}3)--(2bar{1}3)--(213)--(123)--(132)--(231)--
				(23bar{1})--(13bar{2})--(12bar{3})--(bar{1}2bar{3})--(bar{2}1bar{3})--(bar{2}bar{1}bar{3})--(bar{1}bar{2}bar{3})--(1bar{2}bar{3})--(1bar{3}bar{2})--(2bar{3}bar{1})--(2bar{3}1)
				--(1bar{3}2)--(1bar{2}3);
				\draw[solid,line width=0.8pt] (1bar{2}bar{3})--(2bar{1}bar{3})--(21bar{3})--(12bar{3});
				\draw[solid,line width=0.8pt] (2bar{1}bar{3})--(3bar{1}bar{2})--(3bar{2}bar{1})--(3bar{2}1)--(3bar{1}2)--(2bar{1}3);
				\draw[solid,line width=0.8pt] (21bar{3})--(31bar{2})--(32bar{1})--(321)--(312)--(213);
				\draw[solid,line width=0.8pt] (31bar{2})--(3bar{1}bar{2});
				\draw[solid,line width=0.8pt] (32bar{1})--(23bar{1});
				\draw[solid,line width=0.8pt] (321)--(231);
				\draw[solid,line width=0.8pt] (1bar{3}2)--(bar{1}bar{3}2);
				\draw[solid,line width=0.8pt] (1bar{3}bar{2})--(bar{1}bar{3}bar{2});
				\draw[solid,line width=0.8pt] (3bar{2}bar{1})--(2bar{3}bar{1});
				\draw[solid,line width=0.8pt] (3bar{2}1)--(2bar{3}1);
				\draw[solid,line width=0.8pt] (3bar{1}2)--(312);
				\draw[dashed,line width=0.8pt] (bar{2}bar{1}bar{3})--(bar{3}bar{1}bar{2})--(bar{3}bar{2}bar{1})--(bar{3}bar{2}1)--(bar{3}bar{1}2)--(bar{2}bar{1}3)--(bar{1}bar{2}3);
				\draw[dashed,line width=0.8pt] (bar{2}bar{1}3)--(bar{2}13)--(bar{1}23)--(123);
				\draw[dashed,line width=0.8pt] (bar{1}23)--(bar{1}32)--(bar{2}31)--(bar{2}3bar{1})--(bar{1}3bar{2})--(13bar{2});
				\draw[dashed,line width=0.8pt] (bar{1}2bar{3})--(bar{1}3bar{2});
				\draw[dashed,line width=0.8pt] (bar{2}1bar{3})--(bar{3}1bar{2})--(bar{3}2bar{1})--(bar{3}21)--(bar{3}12)--(bar{2}13);
				\draw[dashed,line width=0.8pt] (bar{1}32)--(132);
				\draw[dashed,line width=0.8pt] (bar{2}31)--(bar{3}21);
				\draw[dashed,line width=0.8pt] (bar{2}3bar{1})--(bar{3}2bar{1});
				\draw[dashed,line width=0.8pt] (bar{3}12)--(bar{3}bar{1}2);
				\draw[dashed,line width=0.8pt] (bar{3}1bar{2})--(bar{3}bar{1}bar{2});
				\draw[dashed,line width=0.8pt] (bar{3}bar{2}1)--(bar{2}bar{3}1);
				\draw[dashed,line width=0.8pt] (bar{3}bar{2}bar{1})--(bar{2}bar{3}bar{1});
			\end{tikzpicture}
		}
	\end{minipage}
	\caption{The left weak orders on $S_4$ (left) and $B_3$ (right), where the non-seperable elements are highlighted in light gray, and the elements corresponding under the map $\varphi_3$ are labeled with the same number.\label{leftbruhatorderonB3}}
\end{figure}

We now show that $\varphi_n$ is indeed an isomorphism of posets under the left weak order, as illustrated by Figure \ref{leftbruhatorderonB3}.

Let $w=w_1w_2\cdots w_n$ be a signed permutation of $n$. 
\delete{An integer $i\in[n]$ is called a \emph{negative index} of $w$ if $w_i<0$. A pair $(i,j)$ with $1\leq i<j\leq n$ is called an \emph{inversion} of $w$ if $w_i>w_j$, and it is called a \emph{negative sum pair} of $w$ if $w_i+w_j<0$.} Define
\begin{align*}
	\text{Neg}(w) &= \lbrace i \in [n] \mid w_i < 0 \rbrace,\\
	\text{Inv}(w) &= \lbrace (i, j) \in [n]^2 \mid i < j, w_i > w_j \rbrace,\\
	\text{Nsp}(w) &= \lbrace (i, j) \in [n]^2 \mid i < j, w_i + w_j < 0 \rbrace.
\end{align*}
It is well-known that 
\begin{align}\label{ell(w)b}
	\ell(w)=\#\text{Neg}(w)+\#\text{Inv}(w)+\#\text{Nsp}(w).	
\end{align}
If $w$ is a permutation, then $\text{Neg}(w)$ and $\text{Nsp}(w)$ are empty, so we have $\ell(w)=\#\text{Inv}(w)$ for short.

For the symmetric group $S_n$ and the hyperoctahedral group $B_n$, the left weak order $\leq_L$ can be described combinatorially as follows.
\begin{lemma}[\cite{BB05}, Section 3.1]\label{u<vforS}
	Let $ u, v \in  S_n $. Then $ u \leq_L v $ if and only if \, $ \mathrm{Inv}(u) \subseteq \mathrm{Inv}(v)$.
\end{lemma}

\begin{lemma}[\cite{Yu23}, Theorem 3.1]\label{u<vforB}
	Let $ u, v \in  B_n $. Then $ u \leq_L v $ if and only if $ \mathrm{Inv}(u) \subseteq \mathrm{Inv}(v) $, $ \mathrm{Neg}(u) \subseteq \mathrm{Neg}(v) $, and $\mathrm{Nsp}(u) \subseteq \mathrm{Nsp}(v) $.
\end{lemma}

To obtain the main result, we first develop some properties of the bijection $\varphi_n$.
\begin{lemma}\label{lem:varphinw0sout}
If $w \in K( S_{n+1})$, then $\varphi_n(w_0^{ S_{n+1}}w)=w_0^{ B_n} \varphi_n (w)$.
\end{lemma}

	\begin{proof}
		First assume that $w=w_1\cdots w_{n+1}$ has a global descent at the  position $p$. Then $w=\mathrm{st}(w_1\cdots w_{p})\ominus (w_{p+1}\cdots w_{n+1})$, and hence
		$$w_0^{ S_{n+1}}w=w_0^{ S_{p}}\mathrm{st}(w_1\cdots w_p)\oplus w_0^{ S_{n+1-p}}(w_{p+1}\cdots w_{n+1}).$$ 
		Thus, by Eq.\,\eqref{eq:defnvarphin}, 
		\begin{align*}
		\varphi_n(w_0^{ S_{n+1}}w)
		&=  \varphi_{p-1}\left(w_0^{ S_{p}}\mathrm{st}(w_1\cdots w_p)\right) \oplus  w_0^{ S_{n+1-p}}(w_{p+1}\cdots w_{n+1}).
		\end{align*}
		Clearly by induction on $n$ gives that $\varphi_{p-1}\left(w_0^{ S_{p}}\mathrm{st}(w_1\cdots w_p)\right)=w_0^{ B_{p-1}}\varphi_{p-1}(\mathrm{st}(w_1\cdots w_p))$. Thus,
		\begin{align*}
		\varphi_n(w_0^{ S_{n+1}}w)
		&=w_0^{ B_{p-1}}\varphi_{p-1}(\mathrm{st}(w_1\cdots w_p)) \oplus  w_0^{ S_{n+1-p}}(w_{p+1}\cdots w_{n+1})\\
		&=w_0^{ B_{n}}\left(\varphi_{p-1}(\mathrm{st}(w_1\cdots w_p)) \oplus  w_0^{ B_{n+1-p}}w_0^{ S_{n+1-p}}(w_{p+1}\cdots w_{n+1})\right) \\
		&=w_0^{ B_{n}}\varphi_n(w).
		\end{align*}
		
		Now suppose $p$ is a global ascent of $w$. Then $p$ is a global descent
		of $w_0^{B_n}w$, so that 
		\begin{align*}
			\varphi_n(w)=\varphi_n(w_0^{S_{n+1}}w_0^{S_{n+1}}w)
			=w_0^{B_{n}}\varphi_n(w_0^{S_{n+1}}w),		
		\end{align*}
		from which we get $\varphi_n(w_0^{ S_{n+1}}w)=w_0^{ B_n} \varphi_n (w)$,
		 completing the proof.
	\end{proof}

\begin{lemma}\label{lem:varphinorderpreserves}
	Let $u,v \in K( S_{n+1})$. Then $u \leq_L v$ in $S_{n+1}$ if and only if $\varphi_n(u) \leq_L \varphi_n(v)$ in $B_{n}$.
\end{lemma}
	\begin{proof}
		$(\Rightarrow)$
We proceed by induction on $n$ using recurrence Eq.\,\eqref{eq:defnvarphin}. It is easy to check that the statement is true for $n\leq 1$.
Suppose that $n\geq2$ and the statement holds for all indices less than $n$.
Let $u \leq_L v$ be separable permutations of $ S_{n+1}$. Then, by Lemma \ref{u<vforS}, $\mathrm{Inv}(u)\subseteq \mathrm{Inv}(v)$, which implies 
$\mathrm{GDes}(u)\subseteq \mathrm{GDes}(v)$.

Before we proceed with the proof, we first provide a direct sum decomposition for 
$u$ and $v$.
Let $p\in[n]$ be a global ascent or global descent of $u$. It follows from Eq.\,\eqref{eq:defnvarphin} that 
$\varphi_n(u)=a_u \oplus a'_u$ if $p\in\mathrm{GAsc}(u)$, and 
$\varphi_n(u)=b_u \oplus b'_u$ if $p\in\mathrm{GDes}(u)$, where
\begin{align*}
	a_u&=\varphi_{p-1}(u_1\cdots u_p), &a_u'&=\mathrm{st}(u_{p+1}\cdots u_{n+1}),\\
	b_u&=\varphi_{p-1}(\mathrm{st}(u_1\cdots u_{p})), &b_u'&= w_0^{ B_{n+1-p}} w_0^{ S_{n+1-p}}(u_{p+1}\cdots u_{n+1}).
\end{align*} 
Clearly, $a_u'$ is separable.
By the induction hypothesis, $a_u$ and $b_u$ are both separable signed permutations, while by Lemma \ref{lem:ww-1rl}, $b_u'$ is also separable.
Similarly, for any global ascent or global descent, say $q$, of $v$, we can define separable signed permutations
$a_v,a'_v,b_v,b'_v$ such that
$\varphi_n(v)=a_v \oplus a'_v$ if $q\in\mathrm{GAsc}(v)$, and 
$\varphi_n(v)=b_v \oplus b'_v$ if $q\in\mathrm{GDes}(v)$.

Consider now whether $\mathrm{GDes}(u)$ and $ \mathrm{GDes}(v)$ are empty sets or not.

{\bf Case 1.} $\mathrm{GDes}(u)\neq \emptyset$.
Assume that $p\in\mathrm{GDes}(u)$, and hence $p\in\mathrm{GDes}(v)$. We can take $q=p$ so that
$\varphi_n(u)=b_u \oplus b_u'$, 
$\varphi_n(v)=b_v \oplus b_v'$.
Since $u\leq_L v$, it follows from Lemma \ref{u<vforS} that $\mathrm{st}(u_1\cdots u_{p})\leq_L \mathrm{st}(v_1\cdots v_{p})$, $u_{p+1}\cdots u_{n+1}\leq_L v_{p+1}\cdots v_{n+1}$, and hence $b_u\leq_Lb_v$, $b_u'\leq_L b_v'$. 
Note that all entries of $b_u'$ and $b_v'$ are negative.
Thus by Lemmas \ref{u<vforS} and \ref{u<vforB},
\begin{align*}
	\mathrm{Inv}(\varphi_n(u))&=\mathrm{Inv}(b_u)\cup\mathrm{Inv}(b_u')\cup\{(i,j)\mid 1\leq i\leq p-1<j\leq n\}\\
	&\subseteq \mathrm{Inv}(b_v)\cup\mathrm{Inv}(b_v')\cup\{(i,j)\mid 1\leq i\leq p-1<j\leq n\}\\
	&=\mathrm{Inv}(\varphi_n(v)),
\end{align*}
and similarly, 
\begin{align*}
	\mathrm{Nsp}(\varphi_n(u))&=\mathrm{Nsp}(b_u)\cup\mathrm{Nsp}(b_u')\cup\{(i,j)\mid 1\leq i\leq p-1<j\leq n\}\subseteq  \mathrm{Nsp}(\varphi_n(v)),\\
	\mathrm{Neg}(\varphi_n(u))&=\mathrm{Neg}(b_u)\cup\mathrm{Neg}(b_u')\subseteq  \mathrm{Neg}(\varphi_n(v)).
\end{align*}
This yields $\varphi_n(u)\leq_L\varphi_n(v)$ by Lemma \ref{u<vforB}. 

{\bf Case 2.} $\mathrm{GDes}(u)= \emptyset$, $\mathrm{GDes}(v)\neq\emptyset$.
Then $\mathrm{GAsc}(u)\neq\emptyset$. Suppose that $p\in \mathrm{GAsc}(u)$ and $q\in \mathrm{GDes}(v)$. So
$$
u=u_1\cdots u_p\oplus \mathrm{st}(u_{p+1}\cdots u_{n+1}),\qquad
v=\mathrm{st}(v_1\cdots v_q)\ominus v_{q+1}\cdots v_{n+1},
$$
and $\varphi_n(u)=a_u \oplus a'_u$, $\varphi_n(v)=b_v \oplus b'_v$.
First consider the case where $p\leq q$.
Since $v$ is a separable permutation, it follows that $\mathrm{st}(v_1\cdots v_q)$ and hence $\mathrm{st}(v_1\cdots v_p)\oplus\mathrm{st}(v_{p+1}\cdots v_q)$ avoids $3142$ and $2413$, so they are both separable. Similarly, $u_1\cdots u_p\oplus\mathrm{st}(u_{p+1}\cdots u_q)$ is separable.
Comparing the inversion sets yields that
\begin{align*}
u_1\cdots u_p\oplus\mathrm{st}(u_{p+1}\cdots u_q)\leq_L\mathrm{st}(v_1\cdots v_p)\oplus\mathrm{st}(v_{p+1}\cdots v_q)&\leq_L \mathrm{st}(v_{1}\cdots v_{q}).
\end{align*} 
Together with Eq.\,\eqref{eq:defnvarphin} and the induction hypothesis, applying $\varphi_{q-1}$ gives that
	\begin{align*}
	\varphi_{p-1}(u_1\cdots u_p)\oplus\mathrm{st}(u_{p+1}\cdots u_q)\leq_L\varphi_{p-1}(\mathrm{st}(v_1\cdots v_p))\oplus\mathrm{st}(v_{p+1}\cdots v_q)
	\leq_L\varphi_{q-1} (\mathrm{st}(v_{1}\cdots v_{q})).
	\end{align*} 
It follows from Lemma \ref{u<vforB} that
\begin{align*}
	\varphi_{n}(u)&=\varphi_{p-1}(u_1\cdots u_p)\oplus\mathrm{st}(u_{p+1}\cdots u_{n+1})\\
	&\leq_L \varphi_{p-1}(u_1\cdots u_p)\oplus\left(\mathrm{st}(u_{p+1}\cdots u_q)\ominus\mathrm{st}(u_{q+1}\cdots u_{n+1})\right)\\
	&=\left (\varphi_{p-1}(u_1\cdots u_p)\oplus\mathrm{st}(u_{p+1}\cdots u_q)\right)\ominus\mathrm{st}(u_{q+1}\cdots u_{n+1})\\
	&\leq_L\left(\varphi_{p-1}(\mathrm{st}(v_1\cdots v_p))\oplus\mathrm{st}(v_{p+1}\cdots v_q)\right) \oplus w_0^{ B_{n+1-q}} w_0^{ S_{n+1-q}}(v_{q+1}\cdots v_{n+1})\\
	&\leq_L\varphi_{q-1}(\mathrm{st}(v_1\cdots v_q))\oplus w_0^{ B_{n+1-q}} w_0^{ S_{n+1-q}}(v_{q+1}\cdots v_{n+1})\\
	&=\varphi_{n}(v).
\end{align*} 
Now consider the case where $p>q$.
Completely analogous to the previous proof, we have
\begin{align*}
	\varphi_{p-1}(u_1\cdots u_p)
	&\leq_L \varphi_{p-1}\left(\mathrm{st}(u_1\cdots u_q)\ominus\mathrm{st}(u_{q+1}\cdots u_{p})\right)\\
	&= \varphi_{q-1}(\mathrm{st}(u_1\cdots u_q)\oplus w_0^{ B_{p-q}} w_0^{ S_{p-q}}\mathrm{st}(u_{q+1}\cdots u_p),
\end{align*} 
and hence
\begin{align*}
	\varphi_{n}(u_1\cdots u_{n+1})
	&= \varphi_{p-1}(u_1\cdots u_p)\oplus\mathrm{st}(u_{p+1}\cdots u_{n+1})\\
	&\leq_L\left(\varphi_{q-1}(\mathrm{st}(u_1\cdots u_q)\oplus w_0^{ B_{p-q}} w_0^{ S_{p-q}}\mathrm{st}(u_{q+1}\cdots u_p)\right)\oplus\mathrm{st}(u_{p+1}\cdots u_{n+1})\\
	&=\varphi_{q-1}(\mathrm{st}(u_1\cdots u_q)\oplus \left(w_0^{ B_{p-q}} w_0^{ S_{p-q}}\mathrm{st}(u_{q+1}\cdots u_p)\oplus\mathrm{st}(u_{p+1}\cdots u_{n+1})\right).
\end{align*} 
Since $\mathrm{st}(u_1\cdots u_q)\leq_L\mathrm{st}(v_1\cdots v_q)$ and
\begin{align*}
&\,w_0^{ B_{p-q}} w_0^{ S_{p-q}}\mathrm{st}(u_{q+1}\cdots u_p)\oplus\mathrm{st}(u_{p+1}\cdots u_{n+1})\\
\leq_L&\,w_0^{ B_{p-q}} w_0^{ S_{p-q}}\mathrm{st}(u_{q+1}\cdots u_p)\oplus w_0^{ B_{n+1-p}} w_0^{ S_{n+1-p}}\mathrm{st}(u_{p+1}\cdots u_{n+1})\\
=&\,w_0^{ B_{n+1-q}} w_0^{ S_{n+1-q}}\left(\mathrm{st}(u_{q+1}\cdots u_p)\oplus \mathrm{st}(u_{p+1}\cdots u_{n+1})\right)\\
\leq_L&\,w_0^{ B_{n+1-q}} w_0^{ S_{n+1-q}}\mathrm{st}(v_{q+1}\cdots v_{n+1}),
\end{align*} 
it follows that
\begin{align*}
	\varphi_{n}(u_1\cdots u_{n+1})
	\leq_L\varphi_{q-1}(\mathrm{st}(v_1\cdots v_q))\oplus w_0^{ B_{n+1-q}} w_0^{ S_{n+1-q}}\mathrm{st}(v_{q+1}\cdots v_{n+1})
	=\varphi_{n}(v).
\end{align*} 

{\bf Case 3.} $\mathrm{GDes}(u)= \emptyset$, $\mathrm{GDes}(v)=\emptyset$.
Then both $\mathrm{GAsc}(u)$ and $\mathrm{GAsc}(v)$ are nonempty.
Take $p\in \mathrm{GAsc}(u)$ and $q\in \mathrm{GAsc}(v)$, then $p\in \mathrm{GDes}(w_0^{ S_{n+1}}u)$ and $q\in \mathrm{GDes}(w_0^{ S_{n+1}}v)$.
Since $u\leq_L v$, we have $w_0^{ S_{n+1}}v\leq_L w_0^{ S_{n+1}}u$, and hence $q\in \mathrm{GDes}(w_0^{ S_{n+1}}u)$, so we may suppose $p=q$. It follows from
Lemma \ref{lem:varphinw0sout} that
\begin{align*}
	\varphi_n(u)=w_0^{B_n}\varphi_n( w_0^{ S_{n+1}} u)\leq_Lw_0^{B_n}\varphi_n( w_0^{ S_{n+1}} v)=\varphi_n(v).
\end{align*}
Thus, $\varphi_n$ preserves the left weak order.

$(\Leftarrow)$  If $\varphi_n(u) \leq_L \varphi_n(v)$, then $\mathrm{Inv}(\varphi_n(u)) \subseteq \mathrm{Inv}(\varphi_n(v))$ and $\mathrm{Neg} (\varphi_n(u)) \subseteq \mathrm{Neg} (\varphi_n(v))$. By Lemma \ref{prop:w=st0variw}, $\mathrm{Inv} (u)
=\mathrm{Inv} (\mathrm{st}(0\cdot \varphi_n(u))) 
=\mathrm{Inv} (0\cdot \varphi_n(u))$, 
and similarly $\mathrm{Inv} (v)
=\mathrm{Inv} (0\cdot \varphi_n(v))$.
A direct calculation shows that
		\begin{align*}
			\mathrm{Inv} (u)
			&= \{ (1,j+1) \mid j \in \text{Neg}(\varphi_n(u))\} \cup \{ (i+1, j+1) \mid (i,j) \in \mathrm{Inv} (\varphi_n(u)) \} \\
			&\subseteq \{ (1,j+1) \mid j \in \text{Neg}(\varphi_n(v))\} \cup \{ (i+1, j+1) \mid (i,j) \in \mathrm{Inv} (\varphi_n(v)) \} \\
			&= \mathrm{Inv} (v).
		\end{align*}
		Thus $u \leq_L v$, completing the proof.
	\end{proof}
	
Note that $K(S_{n+1})$ is a subset of $S_{n+1}$, so the left weak order $\leq_L$ on $S_{n+1}$ induces a partial order, still denoted by $\leq_L$, on $K(S_{n+1})$, such that $(K(S_{n+1}),\leq_L)$ is an induced subposet of $(S_{n+1},\leq_L)$. So for any $u,v\in K(S_{n+1})$, $u\leq_L v$ in $K(S_{n+1})$ is equivalent to $u\leq_L v$ in $S_{n+1}$.
Similarly, $(K(B_{n}),\leq_L)$ is an induced subposet of $(B_{n},\leq_L)$. 

\begin{theorem}\label{thm:order-preserving}
	For any nonnegative integer $n$, the map
$\varphi_n$ is a poset isomorphism from $K(S_{n+1})$ to $K(B_n)$ under the left weak order.	
\end{theorem}
\begin{proof}
	By Theorem \ref{thm:bijectionKSn+1KBn}, $\varphi_n$ is a bijection from $K(S_{n+1})$ to $K(B_n)$, and by Lemma \ref{lem:varphinorderpreserves}, $\varphi_n$ and its inverse $\psi_n$ are order-preserving.
	Hence, $\varphi_n$ is a poset isomorphism.	
\end{proof}

\section{Explicit formulas for $F(\Lambda_w^B,q)$ and $F(V_w^B,q)$}\label{sec:formulasforFlambdaVB}

In this section, we derive explicit expressions for the rank generating functions $F(\Lambda_w,q)$ and $F(V_w,q)$ of separable elements $w$ in the Weyl group of type $B$, following the approach of Wei \cite{We12}.
We use the separating tree structure of signed permutations to obtain these formulas. 
For any separable element $w$ in a Weyl group $W$, Gaetz and Gao \cite[Theorem 9]{GG20am}  established formulas for the rank generating functions of the weak order ideals $\Lambda_w$ and $V_w$ in terms of the graph associahedron of the Dynkin diagram of $W$. In contrast, our expressions are derived directly from the one-line notation of signed permutations.

To avoid confusion, in this section, for any $w\in S_n$ we use  $\Lambda_w^ S$ (respectively, $V_w^ S$) to denote the interval $[e,w]$  (respectively, $[w,w_0^{ S_n}]$) in $ S_n$, and similarly $\Lambda_w^ B$ (respectively, $V_w^ B$) for the interval $[e,w]$ (respectively, $[w,w_0^{ B_n}]$) in $ B_n$.

\begin{prop}[\cite{We12}, Propositions 2.6 and 2.7] \label{V1}
Let  $ w = w_1  \ldots w_n \in K( S_n) $.
\begin{enumerate}
\item\label{item:V1-GAscp} If $p\in \mathrm{GAsc}(w)$, then
\begin{align*}
F(\Lambda_w^ S,q)&=F(\Lambda_{w_1\cdots w_p}^ S,q) \cdot F(\Lambda_{\mathrm{st}(w_{p+1}\cdots w_{n})}^ S,q),\\
F(V_w^ S,q)&= \binom{n}{p}F(V_{w_1\cdots w_p}^ S,q) \cdot F(V_{\mathrm{st}(w_{p+1}\cdots w_{n})}^ S,q).
\end{align*}  

\item\label{item:V1-GDesp} If $p\in \mathrm{GDes}(w)$, then
\begin{align*}
F(\Lambda_w^ S,q)&= \binom{n}{p} F(\Lambda_{\mathrm{st}(w_1\cdots w_p)}^ S,q) \cdot F(\Lambda_{w_{p+1}\cdots w_{n}}^ S,q),\\
F(V_w^ S,q)&=F(V_{\mathrm{st}(w_1\cdots w_p)}^ S,q) \cdot F(V_{w_{p+1}\cdots w_{n}}^ S,q).
\end{align*}
\end{enumerate}
\end{prop}

We immediately obtain the following corollary for separable signed permutations with no negative entries.
\begin{coro}\label{coro:F(LambdawB)}
	Let  $w\in K( B_n)$ with $\mathrm{Neg}(w)=\emptyset$.
	Then $		F(\Lambda_w^ B, q) = F(\Lambda_w^ S, q)$.
	More precisely,	
	if  $ w=w_1\cdots w_{p} \oplus \mathrm{st}(w_{p+1}\cdots w_{n})$, then  
	\begin{align*} 
	F(\Lambda_w^ B, q) = F(\Lambda_{w_1\cdots w_{p}}^ B, q) \cdot F(\Lambda_{\mathrm{st}(w_{p+1}\cdots w_{n})}^ B, q);
\end{align*}	
if  $w=\mathrm{st}(w_1\cdots w_{p}) \ominus w_{p+1}\cdots w_{n}$, then   
\begin{align*} 
	F(\Lambda_w^ B, q) = \binom{n}{p}F(\Lambda_{\mathrm{st}(w_1\cdots w_{p})}^ B, q) \cdot F(\Lambda_{w_{p+1}\cdots w_{n}}^ B, q).
 \end{align*}
\end{coro}
\begin{proof}	
	Since $\mathrm{Neg}(w)=\emptyset$, it follows that $w\in  S_n $, so the left weak interval   $[e,w]$ in $ B_n$ coincides with the left weak interval $[e,w]$ in $ S_n$. Thus, \begin{align*}
		F(\Lambda_w^B, q) = F(\Lambda_w^ S, q),
	\end{align*}
	and the desired result follows from Proposition \ref{V1}.
\end{proof}

 \begin{lemma}\label{lem:antiisow_0}
 	For any signed permutation $w\in B_n$,
 	there exists an anti-isomorphism $\mu : [w, w_0^{B_n}] \to [e,w_0^{B_n} w]$ defined by $\mu(\pi) = w_0^{B_n} \pi$ for all $\pi \in [w, w_0^{B_n}]$.
 	Moreover, $\ell(w_0^{B_n} w) = n^2 - \ell(w)$.
 \end{lemma}
 \begin{proof}
 	The existence of the anti-isomorphism $\mu$ follows from \cite[Proposition 3.1.5]{BB05}. Since $\ell(w_0^{B_n})=n^2$, it follows that
 	$\ell(w_0^{B_n} w)=\ell(w_0^{B_n})-\ell(w)= n^2 - \ell(w)$. 
 \end{proof}

The following lemma provides a factorization of $F(\Lambda_w^{B},q)$ and $F(V_w^{B},q)$ for separable signed permutations $w$ with proper subsets of negative entries. 

\begin{lemma}\label{lem:F(LambdawB=acdotb)}
Let $ w = w_1\cdots w_n \in K( B_n) $, where $\mathrm{Neg}(w)$ is a proper subset of $[n]$. 
\begin{enumerate}
	\item\label{lem:item:wn>0Fdeomp} If  $ w_n > 0 $ and  $p = \max \mathrm{Neg}(w) $, then
	\begin{align*}
		F(\Lambda_w^ B, q) = F(\Lambda_{w_1 \cdots w_p}^ B, q) \cdot F(\Lambda_{\mathrm{sts}(w_{p+1}\cdots w_n)}^ B, q).
	\end{align*} 
	\item\label{lem:item:wn<0Fdeomp} If  $ w_n < 0 $ and  $p = \max \mathrm{Neg}(w_0^{ B_n}w)$, then
	\begin{align*}
		F(V_{w}^ B, q) = F(V_{w_1\cdots w_p}^ B, q) \cdot F(V_{\mathrm{sts}(w_{p+1}\cdots w_n)}^ B,q).
	\end{align*}
\end{enumerate}
\end{lemma}
\begin{proof}
\eqref{lem:item:wn>0Fdeomp}  Since $ w_n > 0 $ and $p = \max \mathrm{Neg}(w) $, it follows from Theorem \ref{lem:separablesignedperm}\eqref{item:wsepsignedstsop} that
	$w=w_1\cdots w_p\oplus \mathrm{sts}(w_{p+1}\cdots w_n)$.
	Let $u\in\Lambda_w^{ B}$. Then $u\leq_Lw$, and hence by Lemma \ref{u<vforB},
	$$\mathrm{Inv}(u)\subseteq \mathrm{Inv}(w)=\mathrm{Inv}(w_1\cdots w_p)\cup\mathrm{Inv}(w_{p+1}\cdots w_n),$$
	which implies that $u_i<u_j$ for all $i\in[p]$ and $j\in[p+1,n]$. 
	Since  $p = \max \mathrm{Neg}(w) $, we have $u_j>0$ for all $j\in[p+1,n]$ and hence $$\mathrm{Neg}(u)\subseteq\mathrm{Neg}(w)=\mathrm{Neg}(w_1\cdots w_p).$$
    Comparing the sets of negative sum pairs yields that
	\begin{align*}
		\mathrm{Nsp}(u)	\subseteq \mathrm{Nsp}(w)=\mathrm{Nsp}(w_1\cdots w_p),
	\end{align*}
	so $|u_i|<u_j$ for all $i\in[p]$ and $j\in[p+1,n]$, 	
	 and hence $u=u_1\cdots u_p\oplus \mathrm{sts}(u_{p+1}\cdots u_n)$.
	To simplify the notation, we write $u'=u_1\cdots u_p$ and $u''=\mathrm{sts}(u_{p+1}\cdots u_n)$, and hence $u=u'\oplus u''$, and similarly for
	$w=w'\oplus w''$.
	By Lemma \ref{u<vforB}, $u\in\Lambda_w^{ B}$ if and only if 
	$u'\leq_L w'$ and $u''\leq_L w''$.
Thus, the mapping $u\mapsto (u', u'')$ provides a bijection from
$\Lambda_w^ B$ to $\Lambda_{w'}^ B \times \Lambda_{w''}^ B$. 
Note that $\ell (u)= \ell (u') +\ell (u'')$. Therefore, 
\begin{align*}
F(\Lambda_w^B, q) = \sum_{u' \leq_{L} w'} q^{\ell(u')} \cdot \sum_{u'' \leq_{L} w''} q^{\ell(u'')}=F(\Lambda_{w_1\cdots w_p}^ B, q) \cdot F(\Lambda_{\mathrm{sts}(w_{p+1}\cdots w_n)}^ B, q).
\end{align*}

\eqref{lem:item:wn<0Fdeomp} By Lemma \ref{lem:antiisow_0}, the mapping $w\mapsto w_0^{ B_n}w$ is an antiautomorphism from $V_{w}^ B$ to $\Lambda_{w_0^{ B_n}w}^ B$, so
\begin{align}\label{eq:F(Vwbstep1)}
	F(V_w^ B, q) =q^{\ell(w_0^{ B_n}w)}\cdot F(\Lambda_{w_0^{ B_n}w}^ B, q^{-1}).
\end{align} 
Since $w_0^{ B_n}w=\bar{w_1}\cdots \bar{w_n}$ and its last entry $\bar{w}_n$ is positive, Item\,\eqref{lem:item:wn>0Fdeomp} gives
\begin{align*}
 F(\Lambda_{w_0^{ B_n}w}^ B, q^{-1})
 &=F(\Lambda_{\bar{w_1} \cdots \bar{w_p}}^ B, q^{-1}) \cdot F(\Lambda_{\mathrm{sts}(\bar{w_{p+1}}\cdots \bar{w_n})}^ B, q^{-1})\\
 &=F(\Lambda_{w_0^{{ B_p}}(w_1 \cdots w_p)}^ B, q^{-1}) \cdot F(\Lambda_{w_0^{{ B_{n-p}}}\mathrm{sts}(w_{p+1}\cdots w_n)}^ B, q^{-1}).
\end{align*} 
By Theorem \ref{lem:separablesignedperm}\eqref{item:wsepsignedstsop}, $w=w_1\cdots w_p\oplus \mathrm{sts}(w_{p+1}\cdots w_n)$, where  $w_j<0$ for all $j\in[p+1,n]$.
Using Eq.\,\eqref{ell(w)b}, we obtain that
$$
\ell(w)=\ell(w_1\cdots w_p)+\ell(\mathrm{sts}(w_{p+1}\cdots w_n))+2p(n-p).
$$
Thus, by Lemma \ref{lem:antiisow_0},
\begin{align*}
	\ell(w_0^{ B_p}(w_1 \cdots w_p))+	\ell(w_0^{ B_{n-p}}\mathrm{sts}(w_{p+1}\cdots w_n))
	&=p^2-\ell(w_1 \cdots w_p)+(n-p)^2-\ell(\mathrm{sts}(w_{p+1}\cdots w_n)\\
	&=n^2-\ell(w)\\
	&=\ell(w_0^{ B_n}w),
\end{align*} 
and hence Eq.\,\eqref{eq:F(Vwbstep1)} gives
\begin{align*}
	F(V_w^ B, q) &=\left(q^{\ell(w_0^{ B_p}(w_1 \cdots w_p))}\cdot F(\Lambda_{w_0^{ B_{p}}(w_1 \cdots w_p)}^ B, q^{-1})\right)
	\cdot
	\left(q^{\ell(w_0^{ B_{n-p}}\mathrm{sts}(w_{p+1}\cdots w_n))}\cdot F(\Lambda^B_{w_0^{ B_{n-p}}\mathrm{sts}(w_{p+1}\cdots w_n)}, q^{-1})\right)\\
	&= F(V_{w_1\cdots w_p}^ B, q) \cdot F(V_{\mathrm{sts}(w_{p+1}\cdots w_n)}^ B,q),
\end{align*} 
as required.
\end{proof}

\begin{lemma}\label{lem:F(Bn,q)}
	$F(B_n,q) = { [n]!(n)!}$, where $[n]!=[1][2] \cdots [n]$, $[i]=1+q+q^2 +\cdots +q^{i-1}$, $(n)!=(1)(2)\cdots(n)$, and $(i)=1+q^i$.
\end{lemma}
\begin{proof}
	It follows from \cite[Theorem 3.15]{Hum90} or \cite{Che55} that $F(B_n,q) = \prod_{i=1}^n \frac{q^{2i}-1}{q-1}$. The desired identity follows at once by using induction on $n$.
\end{proof}

The following result was first proved for Weyl groups of type $A$ by Wei \cite{We12}, and then for all Weyl groups by Gaetz and Gao \cite{GG20aam}.

\begin{lemma}[\cite{GG20aam}, Theorem 3.9]\label{F(W,q)qtypeA}
	Let $W$ be a Weyl group, and let $w\in W$ be separable. Then 
	the weak order ideals $V_w$ and $\Lambda_w$ are both rank symmetric and rank unimodal. Moreover, 
	\begin{align*}
		F(W,q)=F(\Lambda_w,q)F(V_w,q).
	\end{align*} 
\end{lemma}

\begin{lemma}\label{thm:mainformula1}
Let $w\in K( B_n)$. Then
\begin{align}
F(\Lambda_w^ B, q) &= F(\Lambda_{\mathrm{st}(w)}^ S, q) \cdot \prod_{i \in \mathrm{Neg}(w)} (i),\label{eq:F(Lambdabq1)}\\
F(V_w^ B, q)&= F(V_{\mathrm{st}(w)}^ S, q) \cdot \prod_{i \notin \mathrm{Neg}(w)} (i).\label{eq:F(Lambdabq2)}
\end{align}
\end{lemma}
\begin{proof}
By Corollary \ref{lem:st(w)inksn}, the standard permutation of any separable signed permutation is a separable permutation. Thus, Lemma \ref{F(W,q)qtypeA} implies
	\begin{align*}
		F(\Lambda_{\mathrm{st}(w)}^ S, q) \cdot F(V_{\mathrm{st}(w)}^ S, q)=[n]!.
	\end{align*}
By Lemmas \ref{lem:F(Bn,q)} and \ref{F(W,q)qtypeA},
\begin{align*}
F(V_w^ B, q) F(\Lambda_w^ B, q)= [n]!(n)!.
\end{align*}
Since $(n)!=\prod\limits_{i \in \mathrm{Neg}(w)} (i)\cdot \prod\limits_{i \notin \mathrm{Neg}(w)} (i)$,  
we see that Eq.\,\eqref{eq:F(Lambdabq1)} follows at once from Eq.\,\eqref{eq:F(Lambdabq2)} and vice versa.

The proof is by induction on $n$, where the case $n = 1$ being trivial.  Assume that $n\geq2$. First consider the case where $w_n>0$.
If  $\mathrm{Neg}(w)=\emptyset$, then, by Corollary \ref{coro:F(LambdawB)},
 \begin{align*}
F(\Lambda_w^ B, q)=F(\Lambda_w^ S, q)  =  F(\Lambda_{\mathrm{st}(w)}^ S, q) \cdot \prod_{i \in \text{Neg}(w)} (i). 
\end{align*}
If $\emptyset\neq \mathrm{Neg}(w)\subsetneqq [n]$, then let $p=\max\mathrm{Neg}(w)$. So $p<n$ and $w_j>0$ for all $j\in[p+1,n]$. It follows from Lemma \ref{lem:F(LambdawB=acdotb)}\eqref{lem:item:wn>0Fdeomp} that
\begin{align*}
	F(\Lambda_w^ B, q)& = F(\Lambda_{w_1 \cdots w_p}^ B, q) \cdot F(\Lambda_{\mathrm{st}(w_{p+1}\cdots w_n)}^ B, q).
\end{align*}
By the induction hypothesis on $n$ we have
\begin{align*}
	F(\Lambda_w^ B, q)& = F(\Lambda_{\mathrm{st}(w_1 \cdots w_p)}^ S, q)\cdot \prod_{i \in \mathrm{Neg}(w_1 \cdots w_p)} (i)\cdot F(\Lambda_{\mathrm{st}(w_{p+1}\cdots w_n)}^ S, q).	
\end{align*}
Since $\mathrm{st}(w)=\mathrm{st}(w_1 \cdots w_p)\oplus\mathrm{st}(w_{p+1}\cdots w_n)$, Proposition \ref{V1}\eqref{item:V1-GAscp} yields that
\begin{align*}
	F(\Lambda_{\mathrm{st}(w_1 \cdots w_p)}^ S, q)\cdot F(\Lambda_{\mathrm{st}(w_{p+1}\cdots w_n)}^ S, q)=F(\Lambda_{\mathrm{st}(w)}^ S, q).	
\end{align*}
But $ \mathrm{Neg}(w_1 \cdots w_p)= \mathrm{Neg}(w)$, so
\begin{align*}
	F(\Lambda_w^ B, q)& = F(\Lambda_{\mathrm{st}(w)}^ S, q)\cdot \prod_{i \in \mathrm{Neg}(w)} (i).	
\end{align*}
Therefore, Eq.\,\eqref{eq:F(Lambdabq1)} and thus Eq.\,\eqref{eq:F(Lambdabq2)} hold when $w_n>0$.

Now consider the case where $w_n<0$.
 If  $\mathrm{Neg}(w)=[n] $, then 
$V_w^ B=[w,w_0^{ B_n}]$ is anti-isomorphic to $[e,w_0^{ B_n}w]=\Lambda_{w_0^{ B_n}w}^ B$ by Lemma \ref{lem:antiisow_0}. Note that $w_0^{ B_n}w$ is a separable permutation, so $\Lambda_{w_0^{ B_n}w}^ B=\Lambda_{w_0^{ B_n}w}^ S$, which is anti-isomorphic to $V_{w_0^{ S_n}w_0^{ B_n}w}^ S$. Thus $V_w^ B$ is isomorphic to $V_{w_0^{ S_n}w_0^{ B_n}w}^ S$.
Since $w_0^{ S_n}w_0^{ B_n}w$ is a permutation,
we have $w_0^{ S_n}w_0^{ B_n}w=\mathrm{st}(w_0^{ S_n}w_0^{ B_n}w)=\mathrm{st}(w)$, from which we immediately deduce that
\begin{align*}
	F(V_w^ B, q) = F(V_{w_0^{ S_n}w_0^{ B_n}w}^ S, q)
	= F(V_{\mathrm{st}(w)}^ S, q)
	= F(V_{\mathrm{st}(w)}^ S, q) \cdot \prod_{i \notin \mathrm{Neg}(w)} (i). 
\end{align*}
If $\mathrm{Neg}(w)$ is a proper subset of $[n]$, then let $p=\max\mathrm{Neg}(w_0^{ B_n}w)$. So $p<n$ and $w_j<0$ for all $j\in[p+1,n]$.
It follows from Lemma \ref{lem:F(LambdawB=acdotb)}\eqref{lem:item:wn<0Fdeomp} that
\begin{align*}
	F(V_{w}^ B, q) = F(V_{w_1\cdots w_p}^ B, q) \cdot F(V_{\mathrm{sts}(w_{p+1}\cdots w_n)}^ B,q).
\end{align*}
By the induction hypothesis and the fact that $w_{j}<0$ for all $j\in[p+1,n]$, 
\begin{align*}
	F(V_{w_1\cdots w_p}^ B, q)
	&=F(V_{\mathrm{st}(w_1\cdots w_p)}^ S, q) \cdot \prod_{i \notin \mathrm{Neg}(w_1\cdots w_p)} (i)
	=F(V_{\mathrm{st}(w_1\cdots w_p)}^ S, q) \cdot \prod_{i \notin \mathrm{Neg}(w)} (i),\\
	F(V_{\mathrm{sts}(w_{p+1}\cdots w_n)}^ B,q)&=F(V_{\mathrm{st}(\mathrm{sts}(w_{p+1}\cdots w_n))}^ S, q) \cdot \prod_{i \notin \mathrm{Neg}(\mathrm{sts}(w_{p+1}\cdots w_n))} (i)
	=F(V_{\mathrm{st}(w_{p+1}\cdots w_n)}^ S, q).
\end{align*}
Thus,
\begin{align*}
	F(V_{w}^ B, q) = F(V_{\mathrm{st}(w_1\cdots w_p)}^ S, q) \cdot F(V_{\mathrm{st}(w_{p+1}\cdots w_n)}^ S, q)\cdot \prod_{i \notin \mathrm{Neg}(w)} (i).
\end{align*}
Note that 
$$
\mathrm{st}(w)=\mathrm{st}(w_1\cdots w_p)\ominus \mathrm{st}(w_{p+1}\cdots w_n).
$$
By Proposition \ref{V1}\eqref{item:V1-GDesp}, we have
\begin{align*}
	F(V_{w}^ B, q) = F(V_{\mathrm{st}(w)}^ S, q) \cdot \prod_{i \notin \mathrm{Neg}(w)} (i).
\end{align*}
Consequently, Eq.\,\eqref{eq:F(Lambdabq2)} and thus Eq.\,\eqref{eq:F(Lambdabq1)} hold when $w_n<0$.
This completes the proof.
 \end{proof}

The separating tree $T_{w}$ associated with a separable signed permutation $w$ is defined inductively as follows.

Let $w=w_1\cdots w_n$ be a separable signed permutation. 
\begin{enumerate}
	\item[(1)] If $n=1$, then $w$ is either $1$ or $\bar{1}$. The trees  $T_{1}$ and $T_{\bar{1}}$ are  trivial, consisting of a single node labeled by $1$ and $\bar{1}$, respectively. 
	
	\item[(2)] If $n\geq2$, then by Theorem \ref{lem:separablesignedperm}, there
	exists a separating index $p\in[n-1]$ such that $w$ admits one of the following decompositions:
	
	\begin{itemize}
		\item  $w=w_1\cdots w_{p}\oplus \mathrm{sts}(w_{p+1}\cdots w_n)$,
		where both $w_1\cdots w_{p}$ and  $\mathrm{sts}(w_{p+1}\cdots w_n)$ are separable, and $w_{p+1},\ldots, w_n$ have the same sign. If $w_p<w_{p+1}$, then the corresponding separating tree is constructed as 
		$$T_w=(T_{w_1\cdots w_{p}},\oplus,T_{\mathrm{sts}(w_{p+1}\cdots w_n)}),$$
		where the root is labeled by $\oplus$, the left subtree is $T_{w_1\cdots w_{p}}$ with leaves $w_1,\ldots, w_{p}$, and the right subtree is $T_{\mathrm{sts}(w_{p+1}\cdots w_n)}$ with leaves $w_{p+1},\ldots, w_n$.
		If $w_p>w_{p+1}$, then the corresponding separating tree $T_w$ is defined analogously as 
		$$T_w=(T_{w_1\cdots w_{p}},\ominus,T_{\mathrm{sts}(w_{p+1}\cdots w_n)})$$
		where the root is labeled by $\ominus$, the left subtree is $T_{w_1\cdots w_{p}}$ with leaves $w_1,\ldots, w_{p}$, and the right subtree is $T_{\mathrm{sts}(w_{p+1}\cdots w_n)}$ with leaves $w_{p+1},\ldots, w_n$.

		\item $w=\mathrm{sts}(w_1\cdots w_{p})\ominus w_{p+1}\cdots w_n$, where both $\mathrm{sts}(w_1\cdots w_{p})$ and $w_{p+1}\cdots w_n$ are separable,  and $w_{p+1},\ldots, w_n$ have the same sign.  If $w_p<w_{p+1}$, then define
		$$T_w=(T_{\mathrm{sts}(w_1\cdots w_{p})},\oplus,T_{w_{p+1}\cdots w_n});$$
		if $w_p>w_{p+1}$, then define
		$$T_w=(T_{\mathrm{sts}(w_1\cdots w_{p})},\ominus,T_{w_{p+1}\cdots w_n}).$$
	\end{itemize}	
\end{enumerate}

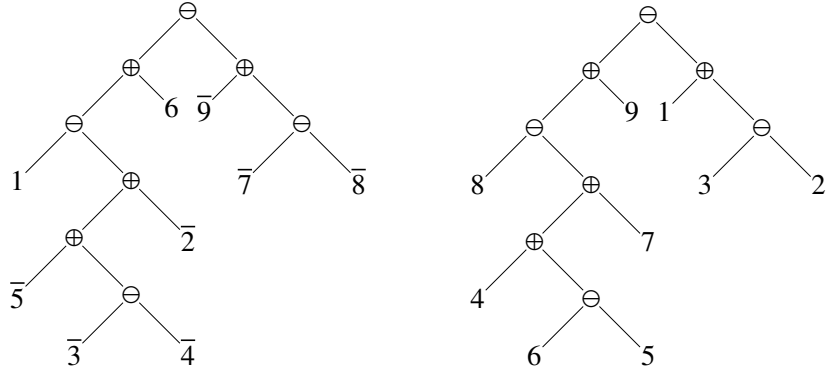
\begin{figure}
	\begin{center}
		\scalebox{0.9}{
	\begin{tikzpicture}[scale=0.6]		
		\node at (-0.2,-0.2) {$1$};
		\draw[-] (0,0) to (1,1);
		\node at(1.2,1.2) {$\ominus$};
		\draw[-] (1.4,1.4) to (2.4,2.4);
		\node at(2.6,2.6) {$\oplus$};
		\draw[-] (2.8,2.8) to (3.8,3.8);
		\node at(4,4) {$\ominus$};
		\draw[-] (4.2,3.8) to (5.2,2.8);
		\node at(5.4,2.6) {$\oplus$};
		\draw[-] (4.6,1.8) to (5.2,2.4);
		\node at(4.4,1.6) {$\bar{9}$};
		\draw[-] (2.8,2.4) to (3.4,1.8);
		\node at(3.6,1.6) {$6$};
		\draw[-] (5.6,2.4) to (6.6,1.4);
		\node at(6.8,1.2) {$\ominus$};
		\draw[-] (7,1) to (8,0);
		\node at(8.2,-0.2) {$\bar{8}$};
		\draw[-] (5.6,0) to (6.6,1);
		\node at(5.4,-0.2) {$\bar{7}$};
		\draw[-] (1.4,1) to (2.4,0);
		\node at(2.6,-0.2) {$\oplus$};
		\draw[-] (2.8,-0.4) to (3.8,-1.4);
		\node at(4,-1.6) {$\bar{2}$};
		\draw[-] (2.4,-0.4) to (1.4,-1.4);
		\node at(1.2,-1.6) {$\oplus$};
		\draw[-] (1,-1.8) to (0,-2.8);
		\node at(-0.2,-3) {$\bar5$};
		\draw[-] (1.4,-1.8) to (2.4,-2.8);
		\node at(2.6,-3) {$\ominus$};
		\draw[-] (2.8,-3.2) to (3.8,-4.2);
		\node at(4,-4.4) {$\bar{4}$};
		\draw[-] (2.4,-3.2) to (1.4,-4.2);
		\node at(1.2,-4.4) {$\bar{3}$};
	\end{tikzpicture}
}
\qquad
\scalebox{0.9}{
		\begin{tikzpicture}[scale=0.6]		
		\node at (-0.2,-0.2) {$8$};
		\draw[-] (0,0) to (1,1);
		\node at(1.2,1.2) {$\ominus$};
		\draw[-] (1.4,1.4) to (2.4,2.4);
		\node at(2.6,2.6) {$\oplus$};
		\draw[-] (2.8,2.8) to (3.8,3.8);
		\node at(4,4) {$\ominus$};
		\draw[-] (4.2,3.8) to (5.2,2.8);
		\node at(5.4,2.6) {$\oplus$};
		\draw[-] (4.6,1.8) to (5.2,2.4);
		\node at(4.4,1.6) {$1$};
		\draw[-] (2.8,2.4) to (3.4,1.8);
		\node at(3.6,1.6) {$9$};
		\draw[-] (5.6,2.4) to (6.6,1.4);
		\node at(6.8,1.2) {$\ominus$};
		\draw[-] (7,1) to (8,0);
		\node at(8.2,-0.2) {$2$};
		\draw[-] (5.6,0) to (6.6,1);
		\node at(5.4,-0.2) {$3$};
		\draw[-] (1.4,1) to (2.4,0);
		\node at(2.6,-0.2) {$\oplus$};
		\draw[-] (2.8,-0.4) to (3.8,-1.4);
		\node at(4,-1.6) {$7$};
		\draw[-] (2.4,-0.4) to (1.4,-1.4);
		\node at(1.2,-1.6) {$\oplus$};
		\draw[-] (1,-1.8) to (0,-2.8);
		\node at(-0.2,-3) {$4$};
		\draw[-] (1.4,-1.8) to (2.4,-2.8);
		\node at(2.6,-3) {$\ominus$};
		\draw[-] (2.8,-3.2) to (3.8,-4.2);
		\node at(4,-4.4) {$5$};
		\draw[-] (2.4,-3.2) to (1.4,-4.2);
		\node at(1.2,-4.4) {$6$};
	\end{tikzpicture}
}
\end{center}
	\caption{Separating trees for $1\bar{5}\,\bar{3}\,\bar{4}\,\bar{2}\,6\,\bar{9}\,\bar{7}\,\bar{8}$ and its standard permutation $846579132$.\label{exmp:septree}}
\end{figure}

A separable signed permutation $w$ may have multiple decompositions of the form $w=u\oplus v$ or $w=u\ominus v$,  leading to different separating trees. To ensure uniqueness, we adopt the decomposition by selecting the largest possible separating index $p$.
This choice guarantees that the underlying binary tree structure coincides with the di-sk tree of the unsigned permutation $\mathrm{st}(w)$, as defined in \cite{FLZ18}.
Furthermore, when $w$ is an ordinary permutation, the tree $T_w$ agrees with the classical separating tree of a permutation \cite{BBL98}.
See Figure \ref{exmp:septree} for the separating trees of $1\bar{5}\,\bar{3}\,\bar{4}\,\bar{2}6\bar{9}\,\bar{7}\,\bar{8}$ and of its standard permutation $846579132$.

Let  $w$ be a separable signed permutation with its corresponding separating tree $T_w$. We define
\begin{align*}
	S^+(w)&=\{V\in T_w\mid  V\ \text{is a positive node whose parent is negative}\},\\
	S^{-}(w)&=\{V\in T_w\mid  V\ \text{is a negative node whose parent is positive}\},
\end{align*}
where the root node of $T_w$ is excluded from both $S^{+}(w)$ and $S^{-}(w)$. For any node $V$ in $T_w$, let 
$N(V)$ denote the number of leaves in the subtree rooted at $V$. In particular, we adopt the empty product convention $\prod_{V\in \emptyset}[N(V)]!=1$.

Given a separable permutation $w$ in $S_n$, the following explicit formulas for the rank generating functions of $\Lambda_{w}^{S}$ and $V_{w}^{S}$ were established by Wei \cite{We12}.

\begin{lemma}[\cite{We12}, Theorem 3.5]\label{F(W,q)qtypeAfom}
Let $w\in K( S_n)$. 
\begin{enumerate}
	 \item If  the root node of $T_w$ is positive, then
\begin{align*}
	F(\Lambda_{w}^{ S},q)=\frac{\prod\limits_{V\in S^{-}(w)}[N(V)]!}{\prod\limits_{V\in S^{+}(w)}[N(V)]!},
	\qquad
	F(V_{w}^{ S},q)=\frac{\prod\limits_{V\in S^{+}(w)}[N(V)]!}{\prod\limits_{V\in S^{-}(w)}[N(V)]!}[n]!.
\end{align*} 

	 \item If the root node of $T_w$ is negative, then
\begin{align*}
	F(\Lambda_{w}^{ S},q)=\frac{\prod\limits_{V\in S^{-}(w)}[N(V)]!}{\prod\limits_{V\in S^{+}(w)}[N(V)]!}[n]!,
	\qquad
	F(V_{w}^{ S},q)=\frac{\prod\limits_{V\in S^{+}(w)}[N(V)]!}{\prod\limits_{V\in S^{-}(w)}[N(V)]!}.
\end{align*} 
\end{enumerate}
\end{lemma}

The main result of this section is as follows.
\begin{theorem}\label{F(W,q)qtypeBfom}
	Let $w\in K( B_n)$. 
	\begin{enumerate}
		\item If $w_1<w_n$, then
		\begin{align*}
			F(\Lambda_{w}^B,q)=\frac{\prod\limits_{V\in S^{-}(w)}[N(V)]!}{\prod\limits_{V\in S^{+}(w)}[N(V)]!}\prod_{i \in \mathrm{Neg}(w)} (i),
			\qquad
			F(V_{w}^B,q)=\frac{\prod\limits_{V\in S^{+}(w)}[N(V)]!}{\prod\limits_{V\in S^{-}(w)}[N(V)]!}[n]!\prod_{i \notin \mathrm{Neg}(w)} (i).
		\end{align*}

		\item\label{item:negat} If $w_1>w_n$, then
		\begin{align*}
			F(\Lambda_{w}^B,q)=\frac{\prod\limits_{V\in S^{-}(w)}[N(V)]!}{\prod\limits_{V\in S^{+}(w)}[N(V)]!}[n]!\prod_{i \in \mathrm{Neg}(w)} (i),
			\qquad
			F(V_{w}^B,q)=\frac{\prod\limits_{V\in S^{+}(w)}[N(V)]!}{\prod\limits_{V\in S^{-}(w)}[N(V)]!}\prod_{i \notin \mathrm{Neg}(w)} (i).
		\end{align*} 
	\end{enumerate}
\end{theorem}
\begin{proof}
According to Lemmas \ref{thm:mainformula1}	and \ref{F(W,q)qtypeAfom}, it suffices to show that $w$ and $\mathrm{st}(w)$ have the same separating tree  when we ignore the labels of leaves. We proceed by induction on $n$, with the case $n=1$ being trivial. Assume that $n\geq2$. 

First consider the case
$w=w_1\cdots w_{p}\oplus \mathrm{sts}(w_{p+1}\cdots w_n)$,
where both $w_1\cdots w_{p}$ and  $\mathrm{sts}(w_{p+1}\cdots w_n)$ are separable, and $w_{p+1},\ldots, w_n$ have the same sign. If $w_1<w_{n}$, 
then $\{w_{p+1},\ldots, w_n\}=[p+1,n]$ and hence $\mathrm{st}(w)=\mathrm{st}(w_1\cdots w_{p})\oplus \mathrm{st}(w_{p+1}\cdots w_n)$. By the induction hypothesis,
$$
T_w=(T_{w_1\cdots w_{p}},\oplus,T_{\mathrm{sts}(w_{p+1}\cdots w_n)})
=(T_{\mathrm{st}(w_1\cdots w_{p})},\oplus,T_{\mathrm{st}(w_{p+1}\cdots w_n)})=T_{\mathrm{st}(w)}.
$$
If $w_1>w_{n}$, then $\{w_{p+1},\ldots, w_n\}=[\bar{n},\overline{p+1}]$  so that
$\mathrm{st}(w)=\mathrm{st}(w_1\cdots w_{p})\ominus \mathrm{st}(w_{p+1}\cdots w_n)$. Consequently,
$$T_w=(T_{w_1\cdots w_{p}},\ominus,T_{\mathrm{sts}(w_{p+1}\cdots w_n)})=(T_{w_1\cdots w_{p}},\ominus,T_{\mathrm{sts}(w_{p+1}\cdots w_n)})=T_{\mathrm{st}(w)}.$$

Now consider the case $w=\mathrm{sts}(w_1\cdots w_{p})\ominus w_{p+1}\cdots w_n$, where both $\mathrm{sts}(w_1\cdots w_{p})$ and $w_{p+1}\cdots w_n$ are separable,  and $w_{p+1},\ldots, w_n$ have the same sign.
By the definition of skew sum, all $w_i$ for $i\in[n]$ have the same sign.
If $w_1<w_{n}$, then 
n.
If $w_1<w_{n}$, then $\{w_{1},\ldots, w_p\}=[\bar{n},\overline{n-p+1}]$,
$\{w_{p+1},\ldots, w_n\}=[\bar{n-p},\overline{1}]$,
we have $\mathrm{st}(w)=\mathrm{st}(w_1\cdots w_{p})\oplus \mathrm{st}(w_{p+1}\cdots w_n)$, so that
$$
T_w=(T_{\mathrm{sts}(w_1\cdots w_{p})},\oplus,T_{w_{p+1}\cdots w_n})
=(T_{\mathrm{st}(w_1\cdots w_{p})},\oplus,T_{\mathrm{st}(w_{p+1}\cdots w_n)})=T_{\mathrm{st}(w)}.
$$
If $w_1>w_{n}$, then $\{w_{1},\ldots, w_p\}=[n-p+1,n]$, $\{w_{p+1},\ldots, w_n\}=[1,n-p]$,
and we have $\mathrm{st}(w)=\mathrm{st}(w_1\cdots w_{p})\ominus \mathrm{st}(w_{p+1}\cdots w_n)$. Hence,
$$
T_w=(T_{\mathrm{sts}(w_1\cdots w_{p})},\ominus,T_{w_{p+1}\cdots w_n})
=(T_{\mathrm{st}(w_1\cdots w_{p})},\ominus,T_{\mathrm{st}(w_{p+1}\cdots w_n)})=T_{\mathrm{st}(w)}.
$$
Therefore, we always have $T_w=T_{\mathrm{st}(w)}$, completing the proof.
	\end{proof}

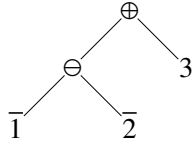
\begin{figure}
	\begin{center}
		\scalebox{0.9}{
			\begin{tikzpicture}[scale=0.6]		
				\node at(2.6,-0.2) {$\oplus$};
				\draw[-] (2.8,-0.4) to (3.8,-1.4);
				\node at(4,-1.6) {$3$};
				\draw[-] (2.4,-0.4) to (1.4,-1.4);
				\node at(1.2,-1.6) {$\ominus$};
				\draw[-] (1,-1.8) to (0,-2.8);
				\node at(-0.2,-3) {$\bar1$};
				\draw[-] (1.4,-1.8) to (2.4,-2.8);
				\node at(2.6,-3) {$\bar2$};
			\end{tikzpicture}
		}
	\end{center}
	\caption{The separating tree $T_{\bar{1}\,\bar{2}3}$.\label{exmp:septreebar123}}
\end{figure}

\begin{exmp} Let $w=1\bar{5}\,\bar{3}\,\bar{4}\,\bar{2}\,6\,\bar{9}\,\bar{7}\,\bar{8}$.
As illustrated by Figure \ref{exmp:septree}, it follows from Theorem \ref{F(W,q)qtypeBfom}\eqref{item:negat} that
	\begin{align*}
		F(\Lambda_{w}^B,q)&=\frac{[2]![2]![5]!}{[3]![4]![6]!}[9]!(2)(3)(4)(5)(7)(8)(9),
		\\
		F(V_{w}^B,q)&=\frac{[3]![4]![6]!}{[2]![2]![5]!}(1)(6).
	\end{align*} 
The separating tree of the separable signed permutation $\bar{1}\,\bar{2}3$ is shown in Figure \ref{exmp:septreebar123}. Consequently, the rank generating functions are computed as follows:
	\begin{align*}
	F(\Lambda_{\bar{1}\,\bar{2}3}^B,q)&=\frac{[2]!}{1}(1)(2)=[2]!(1)(2)=1+2q+2q^2+2q^3+q^4,
	\\
	F(V_{\bar{1}\,\bar{2}3}^B,q)&=\frac{1}{[2]!}[3]!(3)=[3](3)=1+q+q^2+q^3+q^4+q^5.
	\end{align*}	
	These results can be directly verified by the labels of the left weak order shown in Figure \ref{leftbruhatorderonB3}.
\end{exmp}

We conclude this paper with a remark on separable elements of type $D$. The Weyl group of type $D_n$ is isomorphic to the even hyperoctahedral group, consisting of even signed permutations. A natural question arising from our work is to determine the forbidden patterns that characterize separable elements in this group.

\vspace{3mm}
 
{\bf Acknowledgments.}
We are grateful to Shishuo Fu for helpful conversations.
This work was partially supported by the National Natural Science Foundation of China (Grant Nos. 12471020 and 12071377), and the Fundamental Research Funds for the Central Universities (Grant Nos. SWU-XDJH202305 and SWU-KT25015).

\end{document}